\documentclass[11pt]{article}

\usepackage{graphicx}

\newtheorem {theo} {\bf Theorem} [section]
\newtheorem {prop} [theo] {\bf Proposition}
\newtheorem {cory} [theo]{\bf Corollary}
\newtheorem {lem} [theo] {\bf Lemma}

\newtheorem {obs}  [theo]{\bf Remark}
\newtheorem {obss}  [theo]{\bf Remarks}

\newtheorem {afir}  [theo]{\bf Claim}

\newcommand{\CaixaPreta}{\rule{2mm}{2mm}}
\newcommand{\qed}{\hfill\caixapreta \vspace{5mm}}

\newtheorem {REMCURSIVA} [theo] {\bf Remark}

\newcommand{\be}{\begin{eqnarray}}
\newcommand{\ee}{\end{eqnarray}}
\newcommand{\benn}{\begin{eqnarray*}}
\newcommand{\eenn}{\end{eqnarray*}}
\newcommand{\bse}{\begin{equation}}
\newcommand{\ese}{\end{equation}}
\newcommand{\bsenn}{\begin{displaymath}}
\newcommand{\esenn}{\end{displaymath}}

\newcommand{\R}{\mbox{I${\!}$R}}
\newcommand{\hR}{\mbox{I${\!}$R}^{2}_{+}}
\newcommand{\chR}{\overline{\mbox{{I${\!}$R}}}^{2}_{+}}
\newcommand{\Z}{\mbox{Z${\!\!}$Z}}

\newcommand{\rv}{{\bf \mbox {r}}}
\newcommand{\vv}{{\bf \mbox {v}}}
\newcommand{\uu}{{\bf \mbox {u}}}
\newcommand{\xx}{{\bf \mbox {x}}}

\newcommand{\sv}{{\bf \mbox {s}}}

\newcommand{\pp}{{\bf \mbox {p}}}

\title{On the Existence of Periodic Orbits for the Fixed Homogeneous Circle Problem}
 \author{C. Azev\^edo and P. Ontaneda\footnote{The first author was
 partially supported by a CNPq doctorate grant. The second author was
 partially supported by a research grant from CNPq, Brazil.} }
\date{}

\usepackage{wrapfig}
\begin{document}

\maketitle


\begin{abstract} We prove the existence of some types of periodic orbits for
a particle moving in Euclidean three-space under the influence of the gravitational force 
induced by a fixed homogeneous circle. These types include periodic orbits very far and
 very near the homogeneous circle, as well as  eight and spiral periodic orbits.

\end{abstract}

\vspace{0.3cm}

In this paper we use  geometric arguments to demonstrate the existence of
some types of periodic orbits for
the movement of a particle in Euclidean
three-space $\R^{3}$ on which the only acting force is the
gravitational force induced by a fixed homogeneous circle. The study
presented is purely analytical.\\

Interestingly
all we could find in the literature about the fixed homogeneous circle
problem were a few different ways of expressing the potential function. These
expressions appear in classical potential theory books. Among
these expressions are the one expressed in terms of
elliptic integrals of the first kind and the one using the
arithmetic-geometric mean given by Gauss. 
Essentially all expressions of the potential known today had
already appeared in Poincare's Th\'eorie du Potentiel Newtonien
\cite{Po}, published first in 1899. Hence little has been done,
at least in the past century, in the study of this problem.
It is interesting to note that the results proved here use only elementary
geometric constructions (but the technical details are sometimes a
little involving). \\

Before we state our main results we fix some notation. We are
interested in the study of the movement in $\R^{3}$ of a particle
$P$ under the influence of the gravitational force induced by a
fixed homogeneous circle $\cal C$. Denote by ${\mbox{\rv}}=(x,y,z)
\in \R^{3}-{\cal C}$ the position of the particle  $P$ and by
$\dot {\mbox{\rv}}=(\dot x,\dot y,\dot z)$ its velocity. According
to Newton's Law the movement of $P$ obeys the following second
order differential equation:

{\footnotesize\begin{equation} \stackrel{..}{\mbox{\rv}}
\,=\,-\nabla V({\mbox{\rv}}) \label{00}
\end{equation}}

\noindent where  $V$ denotes the potential energy induced by
$\,\cal C$. The expression of $V$ is given by
$V(\mbox{\rv})=-\int_{{\cal
C}}\frac{\lambda\,\,du}{|\!|{\mbox{\rv}}-u|\!|}$,
where  $\lambda$ is the constant mass density of the
circle $\cal C$.
In this introduction we consider the fixed homogeneous circle   {$\,\cal C$}
contained in the $xy$-plane and centered at the origin. Also, by rescaling we can
 consider the circle with radius equal to one (see section 1).
Then the mass of $\cal C$ is given by
$M\, =\, 2\pi\lambda$.

It is intuitively obvious (see section 1.2 of \cite{AO1} for more details)
 that the {\it $z$-axis}, the {\it horizontal plane}
 (i.e. the $xy$-plane, which contains the circle)
and any {\it vertical plane} (i.e any plane that contains the $z$-axis) are invariant
subspaces of our problem. By an invariant subspace $\Lambda$ we mean that
any movement that begins tangentially in $\Lambda$
stays in $\Lambda$ (in the future as well as in the past). Our problem restricted to the
horizontal plane is a central force problem which outside the circle
is given by an attractive force, hence it posses periodic circular orbits. Inside the circle
the force is repulsive and contains no periodic orbits; the dynamics in the horizontal
plane is studied in \cite{AO1}.
\vspace{.1in}

Our first two results show the existence of periodic orbits of the
fixed homogeneous circle problem restricted to a vertical plane.
The first result gives periodic orbits very far and periodic
orbits very close to the circle in any
vertical plane. (In the next Theorem {\it dist} denotes Euclidean distance.)\\

\noindent {\bf Theorem A.} {\it Let $\epsilon>0$ and $\Lambda$ a vertical plane.
Then there exist periodic orbits ${\mbox{\rv}}(t)$, ${\mbox{\sv}}(t)$
in $\Lambda$  of the fixed homogeneous circle problem such that:}
\begin{enumerate}
\item[{i.}] {\it  $|\!|{\mbox{\rv}}(t)|\!|\geq \frac{1}{\epsilon}$, for all $t$ (i.e ${\mbox{\rv}}(t)$
is far from the circle). Moreover, the trace of ${\mbox{\rv}}(t)$ is a simple closed curve,
symmetric  with respect to  the horizontal plane and the $z$-axis,
and encloses the origin.}

\item[{ii.}]  {\it $dist({\mbox{\sv}}(t),{\cal C})\leq \epsilon$, for all $t$
(i.e. $\mbox{\sv}(t)$ is close to the circle).
Moreover, the trace of ${\mbox{\sv}}(t)$ is a simple closed curve, symmetric  with respect
to the horizontal plane,
and encloses the fixed homogeneous circle.}
\end{enumerate}

\begin{figure}[!htb]
 \begin{minipage}[b]{0.43\linewidth}
 \includegraphics[width=\textwidth]{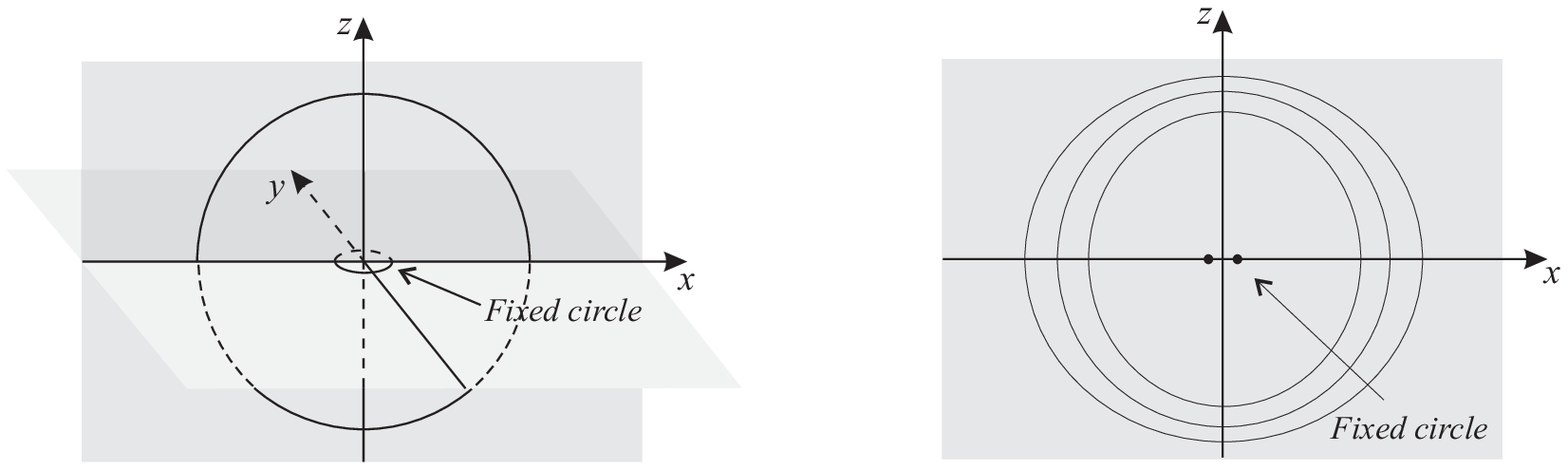}
\caption{\scriptsize{Periodic orbits in the $xz$-plane far from  the fixed homogeneous circle.}}
 \end{minipage}  \hfill
\begin{minipage}[b]{0.43\linewidth}
 \includegraphics[width=\linewidth]{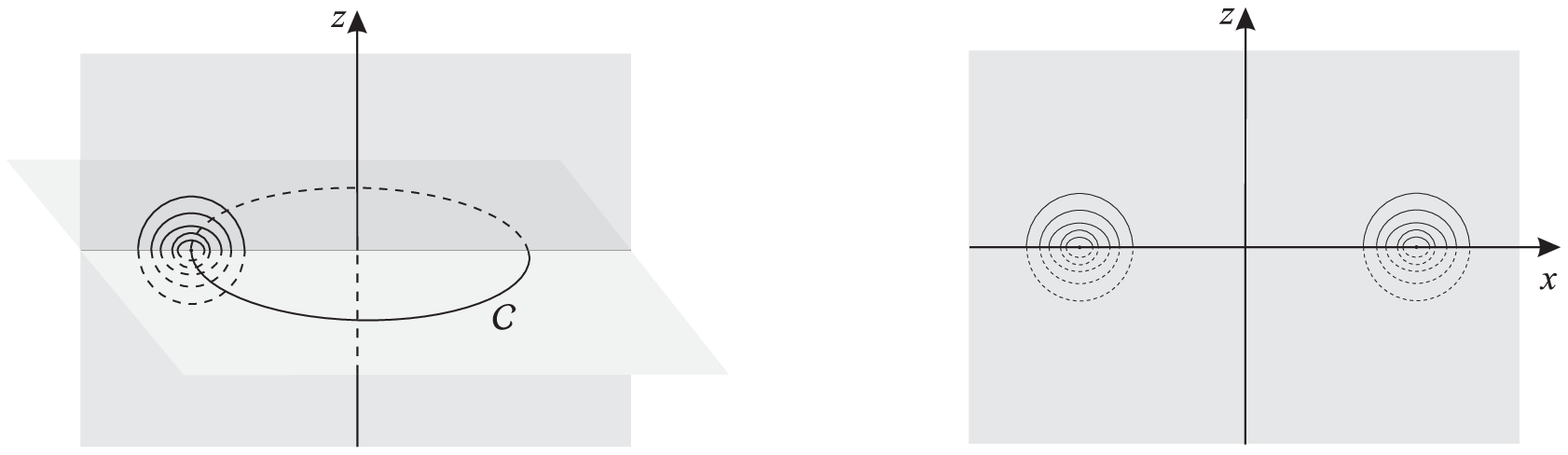}
 \caption{\scriptsize{Periodic orbits in the $xz$-plane near the fixed homogeneous
circle.}}
 \end{minipage}
\end{figure}

Theorem A will follow essentially from a very geometric result about perturbations of certain
central forces \cite {AO}.
Our next result shows the existence of infinitely many ``figure eight" periodic orbits in any vertical plane.
Before we state this result we have to say what we understand for a figure eight orbit.
To fix ideas let the vertical plane $\Lambda$ be the $xz$-plane. Generalizations to any $\Lambda$ are
straightforward. We write $(x,z)$ for the coordinates of a point in the $xz$-plane.
We say that an orbit ${\mbox{\rv}}(t)=(x(t),z(t))$, $t\in [0,\tau ]$, $\tau >0$, of the fixed
homogeneous circle problem
in $\Lambda$ is an {\it essential part of a symmetric figure eight orbit} if
${\mbox{\rv}}(t)$ satisfies:
(1) $z(0)=0,\,\,x(0)>1$, (2) ${\mbox{\rv}}(\tau)=(0,0)$,
(3) $\dot{x}(0)=0$, (4) $z(t)>0,\,\,\,t\in (0,\tau)$.
Note that we do not demand any condition about the angle between $(-1,0)$ and $\dot
{\mbox{\rv}}(\tau).$ This angle can be less than $\frac{\pi}{2}$
as in the figure to the left below, or larger  than
$\frac{\pi}{2}$ as in the figure  to the right below.

\begin{figure}[!htb]
 \centering
 \includegraphics[scale=0.40]{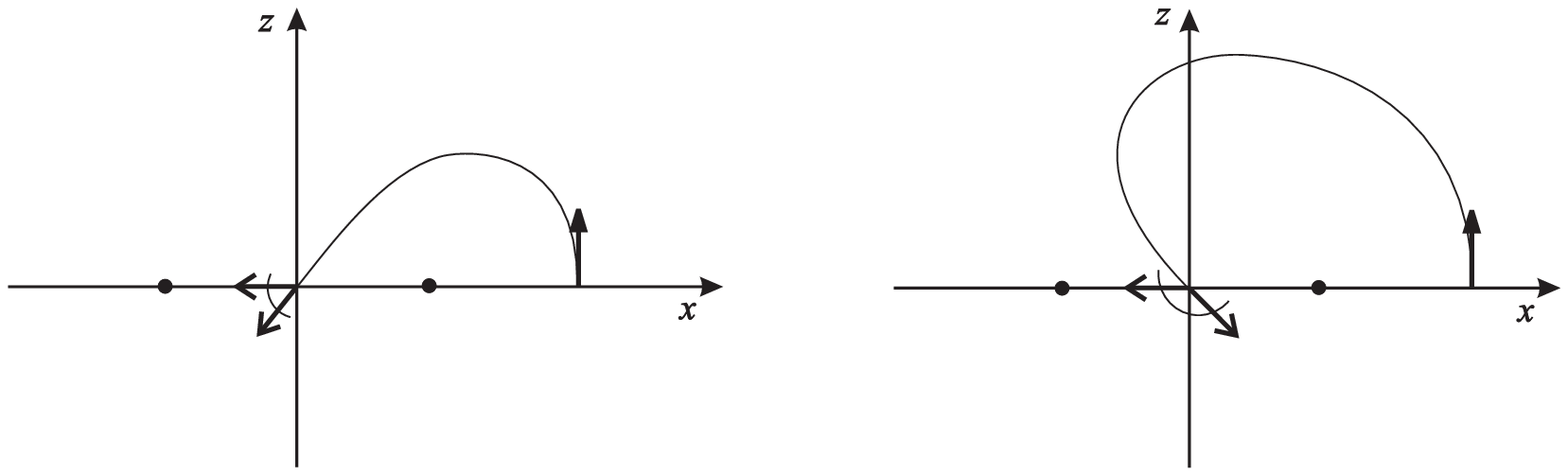}
\caption{\scriptsize{Trace of an essential part of a symmetric figure eight orbit.}}
 \end{figure}

\eject

Now, if $\mbox{\rv}(t)$ is an essential part of a symmetric
figure eight orbit we can use the symmetry of the problem
to define a periodic orbit. Explicitly, using $\mbox{\rv}(t)$
define $\bar {\mbox \rv}:\R\rightarrow\Lambda$
in the following way:

{\footnotesize $$\begin{array}{l}
\bar {\mbox \rv}(t)\,\, =\,\,\left\{
\begin{array}{lll}
{\mbox \rv}(t-4n\tau),& t\in[4n\tau, (4n+1)\tau],&
\\ \varphi_{1}\varphi_{2} {\mbox \rv}((4n+2)\tau-t),
&t\in[(4n+1)\tau,(4n+2)\tau], &
\\
\varphi_{1} {\mbox \rv}(t-(4n+2)\tau),
&t\in[(4n+2)\tau,(4n+3)\tau], &
\\
\varphi_{2} {\mbox \rv}(4n\tau-t), &t\in[(4n-1)\tau,(4n)\tau]. &
\end{array}
\right.
\end{array}$$}

Here $n$ denotes an integer and $\varphi_{1}, \varphi_{2}$ are reflections with
respect to the $z$ and $x$ axes, respectively.
 It is straightforward to verify that $\bar {\mbox \rv}(t)$ is a periodic
orbit of the fixed homogeneous circle problem in $\Lambda$
with  period $4\tau$. Moreover, the trace of $\bar
{\mbox \rv}$ is symmetric with respect to  the $x$ and $z$ axes.
We say that $\bar {\mbox \rv}(t)$ is a {\it symmetric figure eight periodic orbit} in $\Lambda$ and
${\mbox \rv}(t)$ is the essential part of it. \\

\noindent {\bf Theorem B.} {\it In any vertical plane $\Lambda$ there are infinitely many
geometrically distinct symmetric figure eight periodic orbits  of the fixed homogeneous circle
problem, such that their essential parts are imbeddings.}\\

\begin{figure}[!htb]
 \centering
 \includegraphics[scale=0.5]{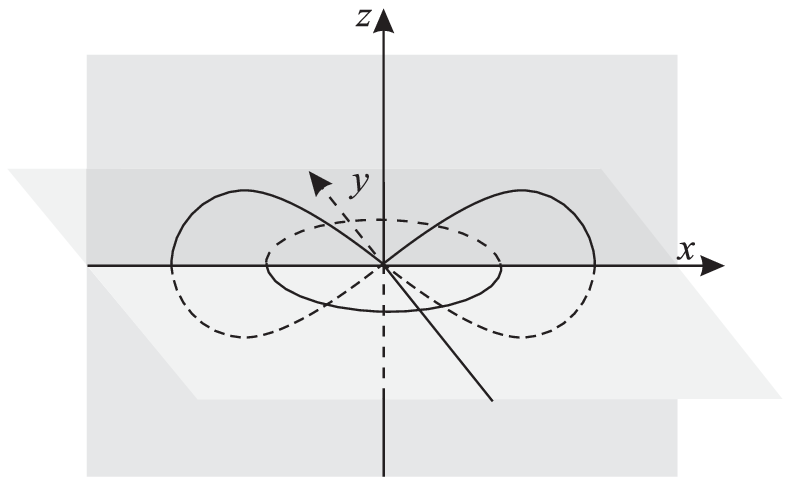}
\caption{\scriptsize{Symmetric figure eight periodic orbit.}}
 \end{figure}

By geometrically distinct orbits we mean orbits whose traces are different.
There is a clear similarity between the the fixed homogeneous circle problem
(restricted to a vertical plane)
and the planar symmetric Euler problem (i.e. the problem of two fixed centers with equal masses).
For the symmetric Euler problem, which is a well-studied
integrable problem, we could not find in the literature any
reference related to the existence of symmetric figure eight orbits. The techniques developed
to prove Theorem B above also work to prove
the existence of symmetric figure eight periodic orbits of the planar symmetric
Euler problem. \\

\noindent {\bf Corollary.} {\it There are infinitely many geometrically distinct
symmetric figure eight periodic orbits of the planar symmetric
Euler problem, such that their essential parts are imbeddings.}\\

Our third result shows the existence of spiral periodic orbits close to the circle.
Before we state this result we have some comments. Since our problem is invariant
by rotations around the $z$-axis, the problem can be reduced in a canonical way
to a problem with two degrees of freedom. In fact,
every orbit ${\mbox \rv}(t)=(x(t),y(t),z(t))$  can be written as
$\mbox{\rv}(t)=(r(t)cos\varphi (t), r(t)sin\varphi(t),z(t))$, i.e the cylindrical
coordinates of ${\mbox \rv}(t)$ are $(r(t),\varphi(t),z(t))$, and
$(r(t),z(t))$ satisfy a system of equations (see (\ref{R4}) in section 5). Moreover, once we know
$r(t)$ we can obtain $\varphi(t)$ by integration. We say that $(r(t),z(t))$ is
the {\it canonical projection orbit} of ${\mbox \rv}(t)$. Geometrically,  $(r(t),z(t))$
is obtained from ${\mbox \rv}(t)$ by a projection ``of a book onto one of its pages":
think of $\R^3$ as book with an infinite number of pages, each page being a half
plane having the $z$-axis as boundary; then identify all pages to one in the
obvious way (like closing the book). Note that under this projection the circle
$\cal C$ projects to a point $x_{\cal C}$ with $rz$ coordinates (1,0).
We say that an orbit is {\it circular} if its canonical projection orbit is
an equilibrium position of the system formed by the first two equations of (\ref{R4}).
 It is proved in \cite{AO1} that all circular orbits lie in
the horizontal plane. Moreover, it is also proved in \cite{AO1} that
 there is a radius $r_0$, $1<r_0<2$ (that does not
depend on the mass) such that a circular orbit in the horizontal
plane with radius $r$ is stable if and only if $r>r_0$. We say
that an orbit (not contained in a vertical plane) is a {\it spiral
orbit} if its canonical projection orbit is periodic. Note that
not every spiral orbit is periodic. Indeed a spiral orbit is
periodic if and only if $\frac{\varphi(\tau)}{2\pi}$ is a rational
number, where $\tau$
is the period of the canonical projection orbit (see Lemma \ref{4.4.2}).\\

\noindent{\bf Theorem C.} {\it For every $\epsilon >0$ and $K\neq 0$ there is a spiral periodic orbit
${\mbox \rv}(t)$ with angular momentum $K$ and such that
$dist({\mbox{\rv}}(t),{\cal C})\leq \epsilon$, for all $t$.
Moreover, its canonical projection orbit is a simple closed curve, symmetric with respect to
the horizontal plane and encloses the
point $x_{\cal C}$.}

\begin{figure}[h]
\centering
\includegraphics[scale=0.35]{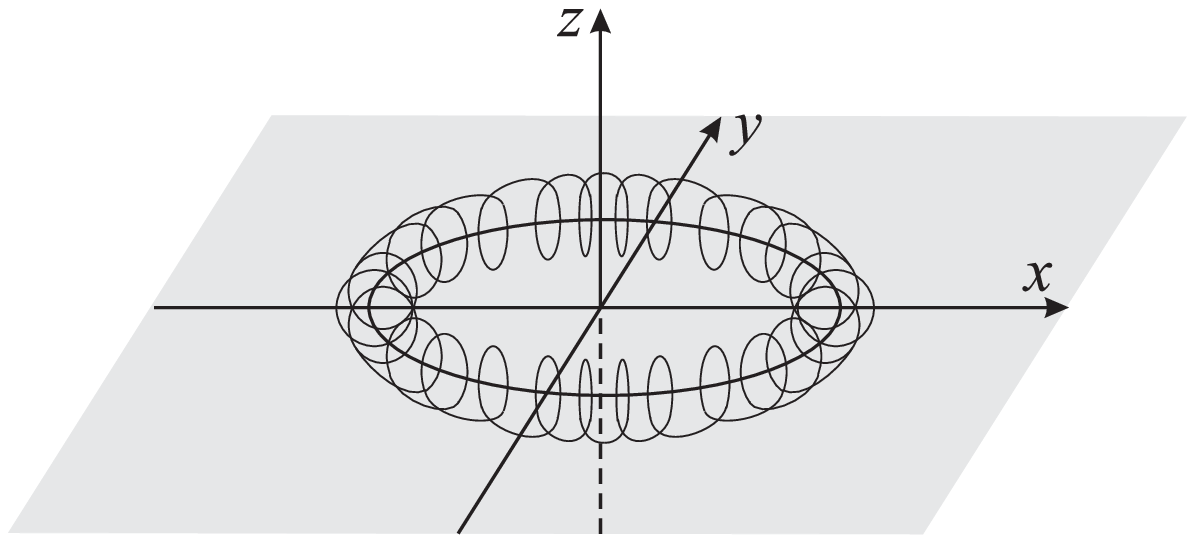}
\caption{\scriptsize{Spiral periodic orbit in $\R^3.$}}
 \end{figure}

Note that the orbits mentioned in Theorem C above lie in a surface of revolution homeomorphic
to a two-torus. It follows easily from the proof of Theorem C that there are also spiral
orbits wich are dense in these tori.
\vspace{0,2cm}

In section 1 we present some preliminary facts. Theorems A, B, C
are proved in sections 2, 3, 5 respectively. In section 4 we show
how use the methods of the proof of Theorem B to prove its
Corollary. Many of our results can be generalized to problems
induced by other symmetric objects. We  give some examples of this
in section 6. Finally, the paper has two appendices. In the first
one we study ``interval pointing forces" in the half plane. The
results in this appendix are used to prove the last part of the
statement of Theorem B: the essential part of the orbits are
imbeddings. These results apply not just to any orbit in a
vertical plane (periodic or not) but to any orbit of a system with
a force that ``points to an interval".
In the second appendix we prove Proposition 3.1.\\

This problem was proposed by H. Cabral, who also
suggested some questions about it.  We are grateful to him. We are
also grateful to T. Stuchi and J. Koiller who suggested some of
the generalizations shown here. Finally, we want to thank M.
Sansuke for drawing most of the figures that appear in this
paper.

\section{Preliminaries.}

We will need to consider circles with variable radius and mass.
Write  $V(\mbox{\rv},\rho,M)$ to denote
the potential induced by $\cal C$ contained in
the $xy$-plane and  centered at the origin, with radius $\rho$ and
mass $M$, and $\nabla
V(\mbox{\rv},\rho,M)$ to denote the gradient (with respect to
$\mbox{\rv}$) of $V(\mbox{\rv},\rho,M)$.

{\lem The   potential $V$ of the  fixed homogeneous
circle problem satisfies the following identities:

\noindent (i) $V(\mbox{\rv},\rho,c\,M)= c
\,V(\mbox{\rv},\rho,M),\,\,$ for  $c\in\R$,

\noindent (ii) $V(c\,\mbox{\rv},c\,\rho,M)= \frac{1}{c} \,
V(\mbox{\rv},\rho,M),\,\,$ for  $c>0$,

\noindent (iii) $\nabla V(c\,\mbox{\rv},c\,\rho,M)=
\frac{1}{c^{2}}\, \nabla V(\mbox{\rv},\rho,M),\,\,$ for
$c>0$,

\noindent (vi) $\nabla V(\mbox{\rv},\rho,c\,M)= c \,\nabla
V(\mbox{\rv},\rho,M),\,\,$ for  $c\in\R$. \label{1.3.1}}

\vspace{0,2cm}

\noindent  {\bf  Proof.} It follows directly from the definition of
$V(\mbox{\rv},\rho,M)=-\frac{M}{2\pi}
\int_{0}^{2\pi}\frac{d\theta}{|\!|\mbox{\rv} -\mbox{$\rho
e^{i\theta}$}|\!|}$, where $e^{i\theta}=(cos\theta,sin\theta,0)$.
\CaixaPreta

{\cory Let  $\rho,\zeta, M, N$ be positive numbers.
If  ${\mbox{\rv}}(t)$ is a solution of
$\,\stackrel{..}{\mbox{\rv}}(t)=-\nabla V({\mbox{\rv}},\rho,M)$
then  ${\mbox{\sv}}(t)=\frac{\zeta}{\rho}\mbox{\rv} \left(
\sqrt{\frac{N\rho^{3}}{M\zeta^{3}}}\,\,t \right)$ is a solution
of  $\,\stackrel{..}{\mbox{\sv}}(t)=-\nabla
V(\mbox{\sv},\zeta,N)$. \label{1.3.2}}

\vspace{0,2cm}

\noindent  {\bf  Proof.} It follows from Lemma \ref{1.3.1}, by a
direct calculation. \CaixaPreta

 {\obs {\rm   Note that if  ${\mbox{\rv}}(t)$ and
${\mbox{\sv}}(t)$ are as above, then they have the same qualitative
properties. For instance, if  ${\mbox{\rv}}(t)$ is periodic (with
period $T$) then ${\mbox{\sv}}(t)$ is also periodic (with
period $\frac{T\sqrt{M\zeta^{3}}}{\sqrt{N\rho^{3}}}$).}
\label{1.3.3}}\\

Note that Lemma \ref{1.3.1} and Corollary  \ref{1.3.2} imply
that in the study the fixed homogeneous circle problem we can assume the
mass and the radius to be equal to one. We will also need the following Lemma.

{\lem Let $\mbox{{\bf \mbox{q}}}\in \R^{n},$
$\Omega\subset\R^{n}$  open and $V:\Omega \rightarrow \R$
 $C^{1}$ in $\Omega$.
If ${\mbox{\rv}}(t)$ is solution of
 $\stackrel{..}{\mbox{\rv}}(t)=-\nabla V({\mbox{\rv}})$ then
${\mbox{\sv}}(t)=\mbox{\rv} (t)+\mbox{{\bf \mbox{q}}}$ is a
solution of  $\,\stackrel{..}{\mbox{\sv}}(t)=-\nabla
W(\mbox{\sv}),$ where $W:\Omega +\mbox{{\bf \mbox{q}}} \rightarrow
\R,$ $\,W(\mbox{\sv}) =V(\mbox{\sv}-\mbox{{\bf \mbox{q}}})$.
Here $\Omega +\mbox{{\bf \mbox{q}}} = \left\{\,\mbox{{\bf
\mbox{\uu}}}+\mbox{\bf \mbox{q}}\,; \,
\mbox{\uu}\in\Omega\,\right\}.$ \label{1.3.4}}

\vspace{0,2cm}

\noindent  {\bf  Proof.} It follows by direct substitution.
\CaixaPreta

\vspace{0,25cm}If $V$ and $W$ are as in the Lemma above, we say
that {\it $W$ is obtained from $V$ by a translation}.

\vspace{0,2cm}

\section{Proof of Theorem A.}

This section has two subsections. In the first one we prove part (i) and in the second one
we prove part (ii).
Without loss of generality we can assume that $\Lambda$ is the $xz$-plane.

\subsection{Periodic solutions   far  from   the fixed homogeneous circle.}

First we show that the  potential of the  fixed
homogeneous circle, with radius $\epsilon$ small, can be regarded
as a perturbation of the potential of the Kepler problem.
Fix $M>0$. Denote by ${\cal C}_{\epsilon}$ the fixed homogeneous
circle contained in the $xy$-plane, centered at the origin,
with radius $\epsilon$ and mass $M$. The potential
$V(\mbox{\rv},\epsilon)$ induced by ${\cal C}_{\epsilon}$ is
given by:

{\footnotesize $$V(\mbox{\rv},\epsilon)
 =-\frac{M}{2\pi}\int_{0}^{2\pi}\frac{d\theta}
{\sqrt{(x-\epsilon\, cos\theta)^{2}+(y-\epsilon\, sin\theta)^{2}+
z^{2}}}$$}

\noindent where $\mbox{\rv}=(x,y,z)$. Note that
$V(\mbox{\rv},\epsilon)$ is defined and analytic in
$\{\,(\mbox{\rv},\epsilon)\,;\, \mbox{\rv}\notin {\cal
C}_{\epsilon}\,\}$. Note also that it makes sense to allow
non-positive values for  $\epsilon$.
In particular $V(\mbox{\rv},0)$ is the potential induced by a point of mass $M$,
located at the origin.
We have that
$V(\mbox{\rv},\epsilon)$ is defined in  $\{ \,(x,y,z,\epsilon)\,;
\,x^2 +y^2 \neq \epsilon^2\,\,\mbox{or}\,\,z\neq 0\,\}$ and
$V$ is even, that is
$V(\mbox{\rv},\epsilon)=V(\mbox{\rv},-\epsilon)$.

{\prop  The potential $V(\mbox{{\mbox{\rv}}},\epsilon)$ induced
by the fixed homogeneous circle ${\cal C}_{\epsilon}$ with radius
$\epsilon$,  is a second order perturbation (with respect to
$\epsilon$) of the Kepler potential, that is,

{\footnotesize $$V(\mbox{{\mbox{\rv}}},\epsilon) = -\frac{M}{\|\mbox{{\mbox{\rv}}}\|}
+\epsilon ^{2}f(\mbox{{\mbox{\rv}}},
\epsilon)$$}
\noindent  where  $ \,f\,$ is analytic in $\{
\,(\mbox{{\mbox{\rv}}},\epsilon) =(x,y,z,\epsilon)\in\R^4\,;
\,0\neq x^2 +y^2 \neq \epsilon^2\,\,\mbox{or}\,\,z\neq 0\,\}$ .
\label{1.6.1}}
\vspace{0,2cm}

\noindent {\bf Proof.} Since
$V(\mbox{{\mbox{\rv}}},\epsilon)=V(\mbox{{\mbox{\rv}}},-\epsilon)$,
a standard analytic continuation argument proves the Proposition.
\CaixaPreta

\vspace{0,2cm} Because $M$ is fixed we have that if $\epsilon\rightarrow 0$ then
$\lambda\rightarrow +\infty.$
Since the problem of the fixed homogeneous  circle, with
fixed mass $M$, can be regarded as a perturbation of Kepler problem we have the following  result:

{\prop
 Let $C$ be a circle in the $xz$-plane,
centered at the origin, and
let   $\,U$ be an open neighborhood of  $C$ of the form
$C\subset U\subset (\R^{2}-\{ (0,0)\})$.
 Then there exists  $\delta_{0}>0$ such  that for each $\epsilon$,
 $0<\epsilon< \delta_{0}$, there exists a periodic solution of
$\stackrel{..}{\mbox{\rv}}\, =-\nabla V({\mbox{\rv}},\epsilon)$
(restricted to the $xz$-plane) in $U$.
 Moreover,  the trace of this solution is a  simple  closed  curve,
 symmetric  with respect to the $x$ and $z$ axes, and encloses the origin.
 \label{4.3.1}}

\vspace{0,2cm}

\noindent {\bf Proof.} The Proposition follows from Proposition \ref{1.6.1} and  Theorem 0.1 of \cite{AO}.
\CaixaPreta\vspace{0,25cm}

 {\obs {\rm   Let ${\mbox{\rv}}(t,\mbox{\xx},\mbox{\vv},\epsilon)$ denote
 a solution of $\stackrel{..}{\mbox{\rv}}\, =-\nabla V({\mbox{\rv}},\epsilon)$
 (restricted to the $xz$-plane) with initial conditions
${\mbox{\xx}}, {\mbox{\vv}}$. Let ${\mbox{\rv}}_0(t)={\mbox{\rv}}(t,\mbox{\xx}_0,\mbox{\vv}_0,0)$
be (circular) solution with trace $C$ and
${\mbox{\rv}}_\epsilon(t)={\mbox{\rv}}(t,\mbox{\xx}_0,\mbox{\vv}_\epsilon,\epsilon)$
a periodic solution given by Proposition \ref{4.3.1}. It follows from Theorem 0.1 of
\cite{AO} that we can choose ${\mbox{\vv}}_\epsilon$ as close
to ${\mbox{\vv}}_0$ as we want.}
\label{2.3}}\\

\noindent {\bf Proof of part (i) of Theorem A.} Let $C$ be a circle as in
Proposition  \ref{4.3.1}, with radius 2 and $U= \{\,
\mbox{\bf
{\mbox\pp}} \in\R^{2}\,;\,{1}<|\!| \mbox{\bf
\mbox{p}}|\!|<{3}\,\}$. By the  Proposition  above there is a
$\delta_{0}>0$ such  that for  each    $\epsilon$, $0<\epsilon<
\delta_{0}$, there exists a solution  ${\mbox{\rv}}_{\epsilon}(t)$ of
 $\stackrel{..}{\mbox{\rv}}\, =-\nabla V({\mbox{\rv}},\epsilon)$
 in $U$ with ${\mbox{\rv}}_{\epsilon}(t)$
periodic  satisfying
$1<|\!|{\mbox{\rv}}_{\epsilon}(t)|\!|<3$, and
whose trace is a simple closed curve,   symmetric with respect to the
$x$ and $z$ axes, and encloses the origin.
By Corollary \ref{1.3.2},
$\mbox{\rv}(t)=\frac{1}{\epsilon}{\mbox{\rv}}_{\epsilon}(\epsilon^{3/2}t)$
 is a solution of
$\stackrel{..}{\mbox{\rv}}\, =-\nabla V({\mbox{\rv}}, 1)$.
Moreover, $\mbox{\rv}(t)$ satisfies
 $\frac{1}{\epsilon}<|\!|{\mbox{\rv}}(t)|\!|<\frac{3}{\epsilon}$,
 for all $t$.
By the properties of ${\mbox{\rv}}_{\epsilon}(t)$ we have that
${\mbox{\rv}}(t)$ is also periodic and  its trace is  a  simple
closed    curve,
 symmetric  with respect to  the $x$ and
 $z$ axes, and encloses the origin.
 \CaixaPreta\vspace{0,2cm}

{\obs {\rm We can assume that the periodic
solutions given by Theorem A (i) intersect the $x$-axis
transversally in exactly  two  points of the form $(x_0,0), (-x_0, 0)$
 (see Remark 2.1 (1) of \cite{AO}).} \label{p1}}\\

Note that if the  period of ${\mbox{\rv}}_{\epsilon}(t)$ is
$\tau$ then  the  period of $\mbox{\sv}(t)$ is  $\tau /
\epsilon^{3/2}$. Hence if $\epsilon$ is small  the  period of
$\mbox{\rv}(t)$ is large. In this way  the periodic orbits
obtained here have large  periods.

\vspace{0,25cm}

\subsection{Periodic solutions  near the fixed homogeneous circle.}

In this subsection instead of fixing the mass $M$ we fix the density
$\lambda >0$. Let   $V(x,y,z;\rho)$  be the  potential induced by the fixed homogeneous
circle, contained in the $xy$-plane, centered at the
origin, with  constant density $\lambda$ and with  radius $\rho$.

\begin{wrapfigure}[5]{o}{6cm}
 \centering
 \includegraphics[width=4cm]{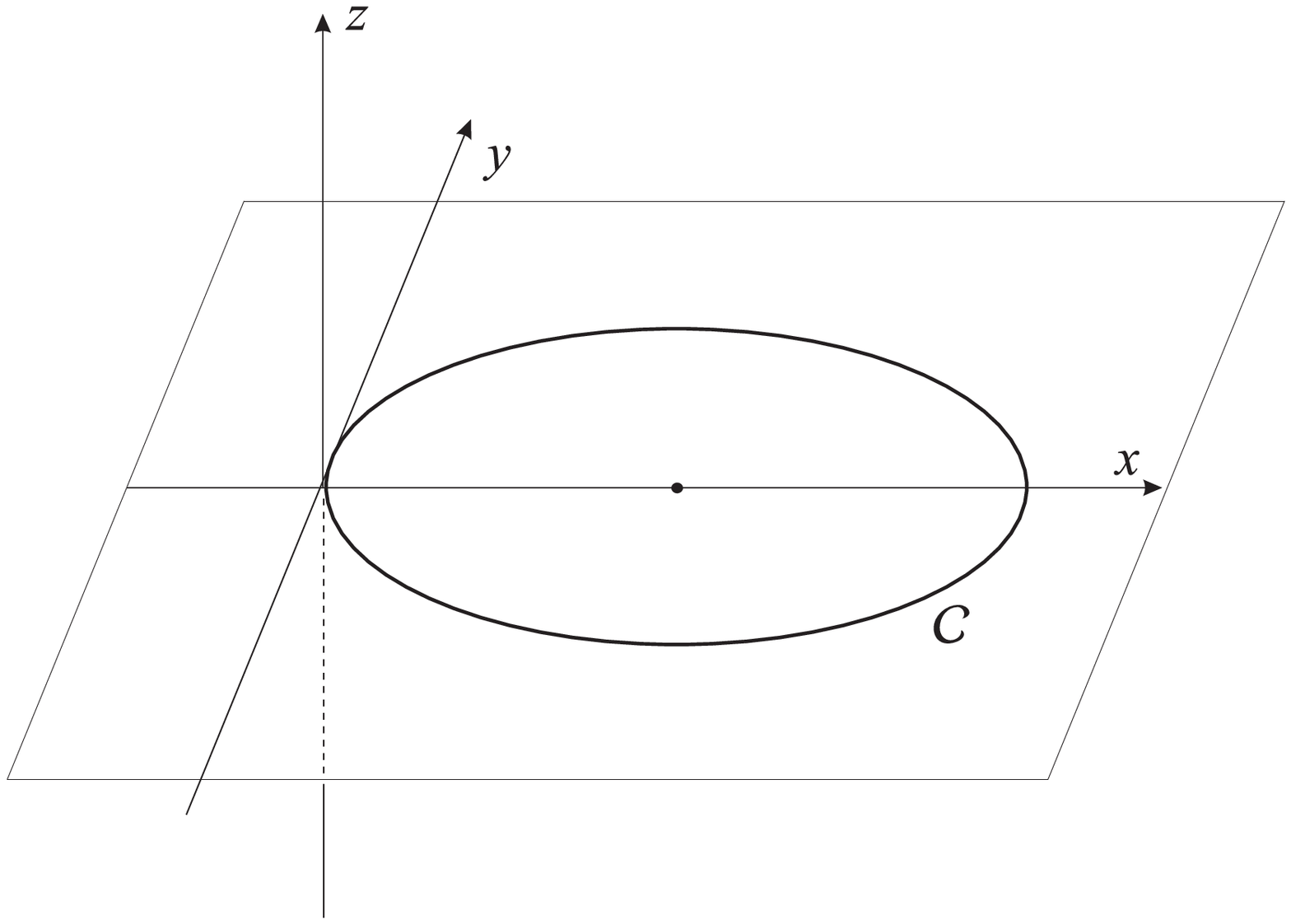}
\caption{\scriptsize{Fixed homogeneous circle centered at $(\frac{1}{\epsilon},0,0)$.}}
\end{wrapfigure}

Now, consider the  fixed homogeneous circle with constant density
$\lambda$, contained in the $xy$-plane and {\it  passing through
the  origin}, with radius $\frac{1}{\epsilon}$ and center
$(\frac{1}{\epsilon},0,0)$, as the figure 2.6. Note that  the
$xz$-plane is a invariant subspace of this new problem.
 Denote  the  potential induced by
this   translated fixed homogeneous circle (and restricted to  the
$xz$-plane)  by $W\left(x,z;\epsilon\right)$. We have that
$W\left(x,z;\epsilon\right)=V\left(
(x,0,z)-(\frac{1}{\epsilon},0,0)\, ;\,\,\,\frac{1}{\epsilon}
\right)$.  That is, $W$ is obtained from   $V$ by a translation.
Note that when  $\epsilon$ tends to zero, the  mass $M=2\pi
\lambda \frac{1}{\epsilon}$ tends to infinity. Let $\nabla W$ be
the gradient of $W$ (with respect to $x,z$). We extend now $\nabla W$ to $\epsilon =0$.
Define $\nabla W(x,z;0) :=64\,\lambda\frac{(x,z)}{x^{2}+z^{2}},$
$(x,z)\neq (0,0).$ Let $A=\{\,(x,z;\epsilon)\,;\, (0,0)\neq (x,z)
\neq (\frac{2}{\epsilon},0)\, \}$. Then $\nabla W(x,z;\epsilon)$ is defined on
$A$.

\vspace{0.2cm}
{\prop  $\nabla W$ is continuous on   $A$.\label{2a}}
\vspace{0.2cm}

The proof of this Proposition appears in \cite{AO1}. Note that
$\nabla W(x,z;0)=64{\lambda}\nabla ln (\sqrt{x^{2}+z^{2}})$ and
recall that $64{\lambda}\,ln (\sqrt{x^{2}+z^{2}})=32{\lambda}ln
(x^{2}+z^{2})$ is the  potential induced by the infinite wire
(with constant density $32\lambda$ and  infinite mass) orthogonal
to $xz$-plane intersecting the $xz$-plane in the   origin. Hence
the problem of the fixed homogeneous circle with large radius,
constant density $\lambda$, can be regarded as a perturbation of
the problem of the infinite homogeneous straight wire with density
$32\lambda$.\vspace{0.2cm}

{\prop
 Let $C$ be a circle  in the  $xz$-plane
with center at origin  and
let   $\,U$ be an open set that contains $C$ of the form
$C\subset U\subset (\R^{2}-\{ (0,0)\})$.
 Then  there is a  $\delta_{0}>0$ such  that for  each    $\epsilon$,
 $0<\epsilon< \delta_{0}$, there exists a periodic solution
in $U$ of $\stackrel{..}{\mbox{\rv}}\, =-\nabla W\left ({\mbox{\rv}},\epsilon\right)$
 (restricted  to the $xz$-plane).
 Moreover, the trace of this solution
is  a simple closed curve,  symmetric  with respect to the x-axis,
and encloses the origin.
\label{4.3.4}}
\vspace{0,2cm}

\noindent {\bf Proof.} It follows from Proposition \ref{2a} and
Theorem 0.2 of \cite{AO}.
\CaixaPreta

\vspace{0,2cm}

{\prop
There is a  $\delta_{0}>0$ such  that for  each    $\epsilon$,
$0<\epsilon< \delta_{0}$, there exists a periodic solution    ${\mbox{\rv}}_{\epsilon}(t)$
of  $ \stackrel{..}{\mbox{\rv}}\, =-\nabla
 V\Bigl({\mbox{\rv}},\frac{1}{\epsilon}\Bigr)$
 (restricted to the $xz$-plane) whose trace is a simple closed curve,
 symmetric  with respect to the x axis,
that encloses the circle with center at the origin and radius $\frac{1}{\epsilon}$.
Moreover, ${\mbox{\rv}}_{\epsilon}(t)$ satisfies
$\frac{1}{3}<dist({\mbox{\rv}}_{\epsilon}(t),{\cal
C}_{\epsilon})<1$.
\label{4.3.5}}
\vspace{0,2cm}

\noindent {\bf Proof.} It follows directly from Lemma \ref{1.3.4} and
from Proposition \ref{4.3.4}, setting $U= \{ \mbox{\bf{
\mbox\pp}}\in\R^{2}\,;\,\frac{1}{3}<|\!|\mbox{\bf
\mbox{p}}|\!|<1\,\}$ and $C$ with radius $\frac{1}{2}$. \CaixaPreta

\vspace{0,2cm}
{\obs {\rm It follows from Theorem 0.2 of \cite{AO} that we can choose
${\mbox{\rv}}_{\epsilon}$ in Proposition \ref{4.3.5} with initial velocity
${\mbox{\vv}}_{\epsilon}=(0, v_\epsilon)$ equal to some fixed velocity
${\mbox{\vv}}_{0}=(0, v_0)$, $v_0 >0$. Also the initial position
${\mbox{\rv}}_{\epsilon}(0)={\mbox{\xx}}_{\epsilon}=(x_\epsilon,0)$
can be such that ${\mbox{\xx}}_{\epsilon}-(\frac{1}{\epsilon},0)$ is
equal to some fixed ${\mbox{\xx}}_{0}=(x_0,0)$, $x_0 >0$ (${\mbox{\xx}}_{0},
{\mbox{\vv}}_{0}$ are the initial conditions of a circular orbit of an unperturbed
problem, see \cite{AO}).} \label{r2.8}}
\vspace{0,4cm}

\noindent {\bf Proof of part (ii) of Theorem A.} Let $\delta_{0}$ and
${\mbox{\rv}}_{\epsilon}(t)$  be as in Proposition  \ref{4.3.5}, with
 $\epsilon< \delta_{0}$.
By Corollary \ref{1.3.2} we have that $\mbox{\sv}(t)=\epsilon\,
{\mbox{\rv}}_{\epsilon}(\frac{t}{\epsilon})$ is a solution of
$\stackrel{..}{\mbox{\rv}}\, =-\nabla V\left ({\mbox{\rv}},1\right)$. Moreover,
$\frac{\epsilon}{3}<dist({\mbox{\sv}}(t),{\cal
{C}})<\epsilon$, and ${\mbox{\sv}}$ satisfies the required properties. \CaixaPreta

\vspace{0,2cm}

{\obs {\rm

\noindent (1) By the  symmetry of the   problem, we have periodic solutions
enclosing the point  $(-1,0)$ and periodic solutions   enclosing the
point  $(1,0)$.

\noindent (2) If the  period  of ${\mbox{\rv}}_{\epsilon}(t)$ is
$\tau$ then the  period  of $\mbox{\sv}(t)$ is $\tau
\epsilon$.
Hence, if $\epsilon$   is  small  the  period of
$\mbox{\sv}(t)$ is small. In this way  the
periodic orbits obtained above have  small period.

 \noindent (3)  We can assume that the periodic
solutions given by Theorem A (ii) intersect the $x$-axis
transversally in exactly  two  points of the form $(x_0,0), (x'_0, 0)$
 with ${x}_0 >1, \,\,0<{x}^{'}_{0}<1$.
 (see Remark 3.2 of \cite{AO}).} \label{r2.9}}

\section{Proof of Theorem B: Figure Eight Orbits.}

Without loss of generality we assume the vertical plane to be the $xz$-plane,
which we identify with $\R^{2}$ with coordinates $(x,z)$. We
define the half vertical plane $\hR =\{ \, (x,z)\in \R^{2}\,;\,
z\,>\, 0\, \}$ and the closed half vertical plane $\chR =\{ \,
(x,z)\in \R^{2}\, ;\, z\, \geq\, 0\, \}$.
In this section we consider the potential
$V({\mbox{\rv}})$ of the circle of radius 1 contained in the $xy$-plane and centered
at the origin. We will use the following notation:  we will write
 $\mbox{\rv}(t, {\mbox\xx}, {\mbox\vv})$ for a solution $\mbox{\rv}(t)$ of
$ \,\stackrel{..}{\mbox{\rv}}\, =-\nabla V({\mbox{\rv}})$
(restricted  to the $xz$-plane) with $\mbox{\rv}(0) ={\mbox\xx},\, \dot{\mbox{\rv}}(0)={\mbox\vv}$.

Note that the potential $V$ has singularities at the points $(-1,0)$
and $(1,0)$, which are the points of intersection of the circle  $\cal C$
with the $xz$-plane.
The system of equations
$ \stackrel{..}{\mbox{\rv}}\, =-\nabla V({\mbox{\rv}})$
 (restricted  to the $xz$-plane) is equivalent to:

{\footnotesize \begin{equation}\left\{ \begin{array}{lll}
\stackrel{..}x&=&-\lambda
\displaystyle\int_{0}^{2\pi}\frac{(x-\,cos\theta)\,
d\theta}{\{\,{(x-\rho\, cos\theta)^{2}+(y-\rho\, sin\theta)^{2}+
z^{2}}\,\} ^{3/2}}
\\ \\
\stackrel{..}z &=&-\lambda \,z \displaystyle\int_{0}^{2\pi}\frac{
d\theta}{\{\,{(x-\rho\, cos\theta)^{2}+(y-\rho\, sin\theta)^{2}+
z^{2}}\,\} ^{3/2}} \label{8.1.4}
\end{array}\right. \end{equation}}

First we mention a property that says that,
under certain conditions, solutions that begin on the $x$-axis
return to it; and they do it with bounded time (depending on the energy).

{\prop For all $\delta <0$ there exists  $T_{\delta}>0$ such
that if $\mbox{\rv}(t)=(x(t),z(t))$ is a maximal solution  of
 $ \,\stackrel{..}{\mbox{\rv}}\, =-\nabla V({\mbox{\rv}})$
(restricted  to the $xz$-plane) with $z(0)=0$, $\dot{z}(0)>0$ and
energy $E(\mbox{\rv})\leq
\delta$, then $\lim_{t\rightarrow t^{-}_{0}}\, z(t)=0$ for some
$t_{0}\in (0,T_{\delta}]$. \label{3.8}}\\

The proof of this Proposition is given in appendix B.
The necessity of writing the limit $\lim_{t\rightarrow
t^{-}_{0}}\, z(t)$ in the statement  above  instead of just
$z(t_0)$ is due to the fact that the solution can approach the circle.
In the proof of Proposition \ref{3.8} we show that
we can choose $T_{\delta} =2\Bigl( 1+\frac{2}{min\{A,1\}} \Bigr),$
where $A=\frac{M}{(1+R_\delta )^{3}}$ and $R_\delta<\infty$ is the radius
of Hill's region of energy $\delta$.
\vspace{0.35cm}

\noindent {\bf Proof of  Theorem B.} From the comments in the
introduction it is enough to prove the existence of infinitely
many essential parts of symmetric figure eight orbits. The fact
that these orbits are imbeddings follows from the results in
appendix A. We will first prove the existence of one essential
part of a symmetric figure eight orbit, and at the end of the
section we indicate how to obtain infinitely many symmetric figure
eight orbits.

We will denote by
 $\mbox\sv _{0}(t)$ a
solution  of $ \,\stackrel{..}{\mbox{\rv}}\, =-\nabla V({\mbox{\rv}})$
 obtained from  Theorem  A (i), that
passes outside (and far from) the fixed homogeneous circle, and
by $\mbox\sv _{1}(t)$  a solution of
$ \,\stackrel{..}{\mbox{\rv}}\, =-\nabla V({\mbox{\rv}})$
 close to the circle, obtained from  Theorem
A (ii), that encloses the point $(1,0)$ (see  figure).
Write $\tilde{\mbox\xx}_{i}=\mbox\sv _{i}(0),$
$\tilde{\mbox\vv}_{i}=\dot{\mbox\sv} _{i}(0),\, i=0,1$.
That is, with the notation introduced above $\mbox\sv _{i}(t)=\mbox\rv
(t,\tilde{\mbox\xx}_{i},\tilde{\mbox\vv}_{i})$.
We will assume the initial conditions $\tilde{\mbox\xx}_{0}=
(\tilde{x}_{0},0)$ and $\tilde{\mbox\xx}_{1}= (\tilde{x}_{1},0)$
in the positive $x$-axis. It is clear that we can also assume that
$\tilde{x}_{0}>\tilde{x}_{1}>1$. Note that the initial velocities
$\tilde{\mbox\vv}_{0}=(0,\tilde{v}_0)$ and
$\tilde{\mbox\vv}_{1}=(0,\tilde{v}_1)$ are orthogonal to the
$x$-axis, since these solutions are symmetric. We can assume that
they have the same direction as the positive $z$-axis, that is
$\tilde{v}_0>0,\,\,\tilde{v}_1 >0$.
By Remark   \ref{p1}, we can choose ${\mbox\sv}_0$ such  that it
intersects transversally the $x$-axis in exactly  two  points
$(\tilde{x}_0, 0)$ and $(-\tilde{x}_0,0)$.
\begin{wrapfigure}[2]{o}{4cm}
 \centering
 \includegraphics[width=3.3cm]{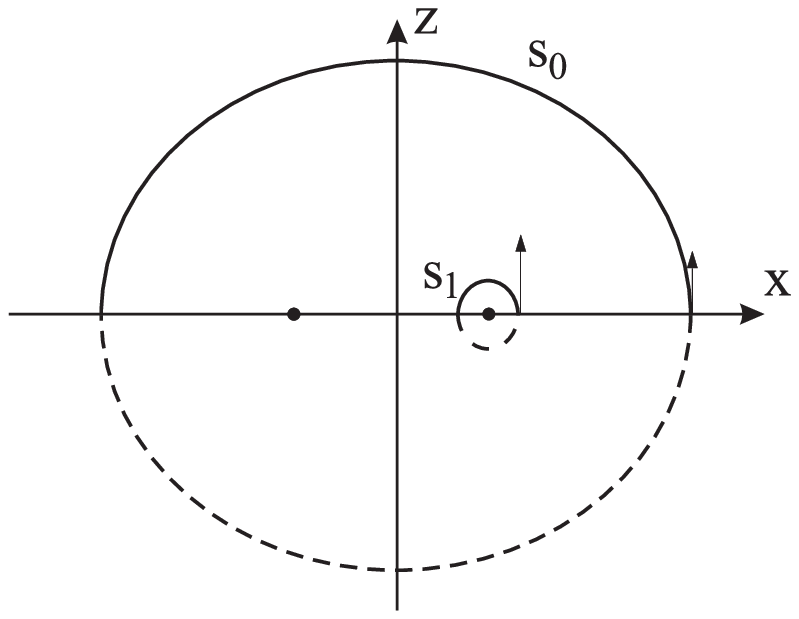}
 \end{wrapfigure}
Also, by Remark   \ref{r2.9} (3), we can choose ${\mbox\sv}_1$ such
that it intersects transversally the $x$-axis   in exactly  two
points $(\tilde{x}_1, 0)$ and $(\tilde{x}^{'}_{1},0),$ with
$0<\tilde{x}^{'}_{1}<1$.

\vspace{0,4cm}

Before we continue with the proof of the Theorem consider following
claims, which we will prove later. Let $h_i =E(\mbox\sv _{i}(t))$ denote
the energy of $\mbox\sv _{i}(t)$, $i=0,1.$ We denote
by $[\tilde{\mbox{\xx}}_{1},\tilde{\mbox{\xx}}_{0}]=\{\,(x,0)\,;\,\tilde{x}_{1}\leq
x\leq \tilde{x}_{0}\,\}$ and
$\,[\tilde{\mbox\vv}_{0},\tilde{\mbox\vv}_{1}]=\{\,(0,v)\,;\,
\tilde{v}_{0}\leq v\leq \tilde{v}_{1}\,\}$.
Let ${\cal A}=
\Bigl([\tilde{\mbox\xx}_{1},\,\tilde{\mbox\xx}_{0}]\times\,\{
\tilde{\mbox\vv}_{0}\}\Bigr)\, \bigcup\,\,
\Bigl(\{\tilde{\mbox\xx}_{1}\}\times\,
[\tilde{\mbox\vv}_{0},\tilde{\mbox\vv}_{1}]\Bigr)$.

{\afir We can choose  $\mbox\sv_{0}(t)$ and $\mbox\sv _{1}(t)$,
 such that
$\tilde{v}_{0}\,< \,\tilde{v}_{1}$ and $h_1<h_0<0$.
\label{1.8}}

{\afir  For  every
$(\mbox\xx,\mbox\vv)\in\cal A$, the energy of $\,\mbox\rv (t)= \mbox{\rv}(t,\mbox\xx,\mbox\vv)$
 is less or equal $h_0$. \label{2.8}}

\vspace{0,3cm}

\begin{wrapfigure}{o}{4cm}
 \centering
 \includegraphics[width=3.5cm]{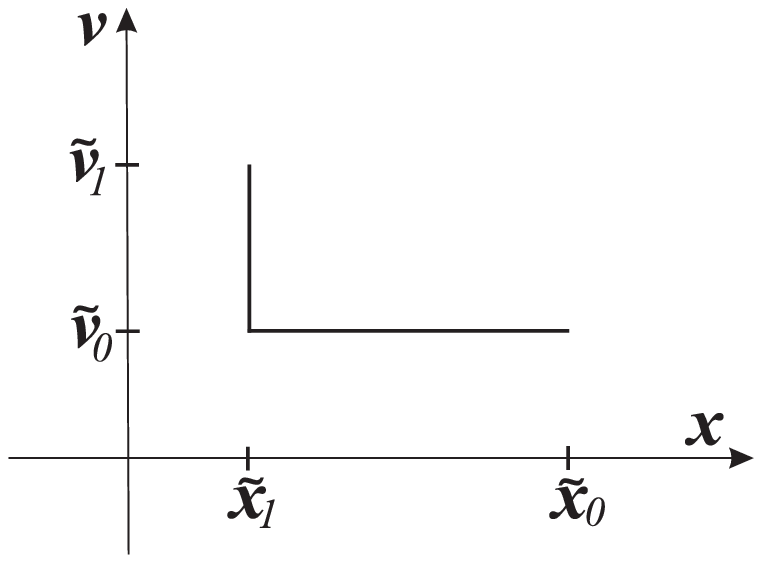}
 \caption{\scriptsize{The set $\{\,(x,v)\,;\, \Bigl(
(x,0),(0,v)\Bigr)\in {\cal A} \,\}$ in the $xv$-plane.}}
 \end{wrapfigure}

Before we state the last claim, we need some comments.
For  $(\mbox\xx,\mbox\vv)\in\cal A$, let   ${\tilde t}
(\mbox\xx,\mbox\vv)= inf\{\,t>0\,;\,\dot{z}(t)=0\,\},$ where
$(x(t),z(t))=\mbox{\rv}(t,\mbox\xx,\mbox\vv)$.
 By Proposition \ref{3.8} and Claim \ref{2.8},
$\tilde{t}(\mbox\xx,\mbox\vv)$ exists and
$\tilde{t}(\mbox\xx,\mbox\vv)\leq T_{h_0},$ for all
$(\mbox\xx,\mbox\vv)\in \cal A$.
By continuity $\dot{z}(\tilde{t}(\mbox\xx,\mbox\vv))=0$.
Since  $\dot{z}(0)=v>0$,
$ \mbox\vv =(0,v),$ we have that
 $\tilde{t}(\mbox\xx,\mbox\vv)>0$.  Hence
$\dot{z}(t) >0$, for all  $t\in
[0,\tilde{t}(\mbox\xx,\mbox\vv))$, and
${z}(t)$ is increasing on
$[0,\tilde{t}(\mbox\xx,\mbox\vv))$.
Define ${\tilde z}(\mbox\xx,\mbox\vv)= {z}({\tilde t}
(\mbox\xx,\mbox\vv))$, with $(\mbox\xx,\mbox\vv)\in\cal A$ and
$(x(t),z(t))=\mbox{\rv}(t,\mbox\xx,\mbox\vv)$.
Since ${\tilde z}(\mbox\xx,\mbox\vv)>0$, the second equation of
(\ref{8.1.4}) implies that ${\tilde z}(\mbox\xx,\mbox\vv)$
is a local maximum of $z(t)$.
\vspace{0,2cm}

The following claim says that there exists a hight $\xi>0$ such
that all solutions  $\mbox{\rv}(t,\mbox\xx,\mbox\vv)$,
with $(\mbox\xx,\mbox\vv)\in\cal A$, ``pass'' this
minimum height.
{\afir There is $\xi >0$ such  that $\,\,{\tilde
z}(\mbox\xx,\mbox\vv)\geq \xi$ for all
$(\mbox\xx,\mbox\vv)\in\cal A$.\label{4.8.8}}

\vspace{0,3cm}

Assuming the claims and  Proposition \ref{3.8} we prove
the Theorem.
Let  $n$ be such  that $\frac{1}{n}<\xi$ where $\xi$ is as  in
Claim  \ref{4.8.8}, and let   $E_{n}$ be the line
$E_{n}=\,\{\,(x,\frac{1}{n})\,;\,x\in\R\,\}$.
 Since $z(0)=0$ we have that
$0=z(0)<\frac{1}{n}<\tilde{z}(\mbox\xx,\mbox\vv)=z(\tilde{t}(\mbox\xx,\mbox\vv)).$
 By  the intermediate value theorem
we have that  there exists  $t_1 =t_1 (\mbox\xx,\mbox\vv)$
such  that
(1)  $0<t_1 (\mbox\xx,\mbox\vv) <\tilde{t}(\mbox\xx,\mbox\vv)$ and
(2) $z(t_1 (\mbox\xx,\mbox\vv))=\frac{1}{n},$ that is
$\mbox{\rv}(t_1 (\mbox\xx,\mbox\vv),\mbox\xx,\mbox\vv)\in E_n.$

Since $\dot{z}(t) >0$ for  $t\in [0,\tilde{t}(\mbox\xx,\mbox\vv))$,
we have  $\dot{z}(t_1
(\mbox\xx,\mbox\vv))>0$. It follows  that $\dot{\mbox{\rv}}(t_1)$
 is   not  horizontal, that is, the intersection of $\mbox{\rv}$ with $E_n$
at $t_1 (\mbox\xx,\mbox\vv)$ is transversal.
For  $(\mbox\xx,\mbox\vv)\in \cal A$ let   $t_2 =t_2
(\mbox\xx,\mbox\vv)=\,inf\{\,t>\tilde{t}(\mbox\xx,\mbox\vv)\,;\,\mbox{\rv}(t,
\mbox\xx,\mbox\vv)\in E_n\,\}.$ Note that   $t_2$ exists
by the intermediate value theorem and Proposition \ref{3.8}. Moreover,
 we have $t_2 (\mbox\xx,\mbox\vv) <T_{h_0}$ and $
z(t)>0,\,\,t\in(0,t_2 (\mbox\xx,\mbox\vv)\,]$.

We show now that $\dot{z}(t_2 (\mbox\xx,\mbox\vv))\neq 0$. If
$\dot{z}(t_2 (\mbox\xx,\mbox\vv))= 0$
 then $t_2 (\mbox\xx,\mbox\vv)$ is local maximum of $z$ (because, since
$z(t_2 (\mbox\xx,\mbox\vv))=1/n >0,\,\, \stackrel{..}z (t_2
(\mbox\xx,\mbox\vv))<0$ by the second equation of (\ref{8.1.4})). It follows  that $z$  is
increasing for $t<t_2 (\mbox\xx,\mbox\vv),$ $t$ near $t_2
(\mbox\xx,\mbox\vv)$.   Since $z(\tilde{t})\geq \xi  >1/n,\,\,$
$\tilde{t}=\tilde{t}(\mbox\xx,\mbox\vv),$ we have that  there exists  $t',
\,$$\tilde{t}<t'<t_2 (\mbox\xx,\mbox\vv)$ with $z(t')=1/n$, which
contradicts the definition of $t_2 (\mbox\xx,\mbox\vv)$. This shows
that $\dot{z}(t_2 (\mbox\xx,\mbox\vv))\neq 0$. It follows that
$\dot{\mbox{\rv}} (t_2 (\mbox\xx,\mbox\vv))$ is not horizontal
thus the intersection  of $\mbox{\rv}$ with $E_n$ at $t=t_2
(\mbox\xx,\mbox\vv) $ is transversal.
Hence, we have functions $t_1, t_2 :{\cal A}\rightarrow \R$
with

(1) $0<t_1 (\mbox\xx,\mbox\vv) <\tilde{t} (\mbox\xx,\mbox\vv) <t_2
(\mbox\xx,\mbox\vv)  <T_{h_0}$,

(2) $\mbox{\rv}(t_i(\mbox\xx,\mbox\vv), \mbox\xx,\mbox\vv)\in
E_n$  and
$\mbox{\rv}(t, \mbox\xx,\mbox\vv)$ intersects $E_n$ transversally
 at $t=t_1, t_2,$

(3)   $t_1,  t_2$ are the first two  times where $\mbox{\rv}(t,
\mbox\xx,\mbox\vv)$ intersects (transversally) $E_n$, that is,
$\mbox{\rv}(t, \mbox\xx,\mbox\vv)\notin E_n,$
 for all $t\in (0,t_2)$, $\,t\neq t_1$.

\vspace{0,23cm}

Note that $t_1$ and $t_2$ depend on $n$.
We show now that $t_1$ and $t_2$ are continuous on
$\cal A$, for  each fixed $n$. Let $({\mbox\xx}_k,
{\mbox\vv}_k)\rightarrow (\mbox\xx,\mbox\vv)$ be a convergent sequence in $\cal A$.
Choose $a,b>0$ such that $t_1
(\mbox\xx,\mbox\vv)  <a<t_2 (\mbox\xx,\mbox\vv)<b$. Then
$\mbox{\rv}(t,\mbox\xx,\mbox\vv),$ $\,a\leq t\leq b,$ intersects
$E_n$ transversally in a single point  (certainly this  point
is $\mbox{\rv}(t_2 (\mbox\xx,\mbox\vv), \mbox\xx,\mbox\vv)$).
Also $\mbox{\rv}(t,\mbox\xx,\mbox\vv),$
$\,0\leq t\leq a,$ intersects   transversally   $E_n$ in a single
point   (this  point is $\mbox{\rv}(t_1 (\mbox\xx,\mbox\vv),
\mbox\xx,\mbox\vv)$). Consequently
 for $k$ sufficiently  large,
$\mbox{\rv}(t,{\mbox\xx}_k,{\mbox\vv}_k),$  $0\leq t\leq a$,
and $\mbox{\rv}(t,{\mbox\xx}_k,{\mbox\vv}_k),$   $a\leq t\leq
b$, also intersect   $E_n$ transversally in a single point
(for this we can use, for example, Proposition 1.4 of \cite{AO}).
It follows  that $\mbox{\rv}(t,{\mbox\xx}_k,{\mbox\vv}_k),$  $0\leq
t\leq b$, intersects   $E_n$ transversally   in two  points and
these intersections happen at times $t_1
({\mbox\xx}_k,{\mbox\vv}_k)$ and $t_2 ({\mbox\xx}_k,{\mbox\vv}_k)$.
By the continuous dependence of solutions of O.D.E.
(e.g. Proposition 1.5 of \cite{AO}) we have
that $\lim_{k\rightarrow +\infty} t({\mbox\xx}_k,{\mbox\vv}_k)=
t(\mbox\xx,\mbox\vv),$ where $t({\mbox\xx}',{\mbox\vv}'),$
$({\mbox\xx}',{\mbox\vv}')$ close to $(\mbox\xx,\mbox\vv)$,
 is the time in which $\mbox{\rv}(t,{\mbox\xx}',{\mbox\vv}'),$
$0\leq t\leq a$, intersects   $E_n$ transversally.
By the uniqueness
of the intersections, $t({\mbox\xx}_k,{\mbox\vv}_k)= t_1
(\mbox\xx_k,\mbox\vv_k),$ which implies that $t_1$ is continuous.
In the same way we have
that $\lim_{k\rightarrow +\infty} t({\mbox\rv}_k
(a),\dot{\mbox\rv}_k (a))= t({\mbox\rv}(a),\dot{\mbox\rv}(a)),$
where ${\mbox\rv}_k (a)=\mbox{\rv}(a,{\mbox\xx}_k,{\mbox\vv}_k) $ and
  $t({\mbox\rv}_k (a),\dot{\mbox\rv}_k (a))$ is the
time  in which $\mbox{\rv}(t,{\mbox\xx}_k,{\mbox\vv}_k),$ $a\leq
t\leq b$, intersects   $E_n$ transversally. By the uniqueness
of the intersections, $t({\mbox\rv}_k (a),\dot{\mbox\rv}_k (a))= t_2
({\mbox\xx}_k,{\mbox\vv}_k)$ and $t({\mbox\rv}(a),
\dot{\mbox\rv}(a))=t_2 (\mbox\xx,\mbox\vv)$. Then    $t_2
({\mbox\xx}_k,{\mbox\vv}_k)\rightarrow t_2
({\mbox\xx},{\mbox\vv}),$  which implies that $t_2$ is continuous.

Now, define  $f:{\cal A}\rightarrow E_{n}$ by
$f({\mbox\xx},{\mbox\vv})={\mbox{\rv}}(t_{2}({\mbox\xx},{\mbox\vv}),
{\mbox\xx},{\mbox\vv})$.
Since $\mbox{\rv}$ and $t_{2}$ are continuous,  $f$  is
continuous. Recall that  $\mbox{\sv}_0 (t)$ and $\mbox{\sv}_1 (t)$
intersect   transversally the $x$-axis  exactly in two  points. Then, for $n$
sufficiently  large,  $\mbox{\sv}_0 (t)$ and $\mbox{\sv}_1 (t)$
also intersect   $E_n$ transversally in exactly two
points, and we have that
$f(\tilde{\mbox\xx}_{0},\tilde{\mbox\vv}_{0})=({x}',\frac{1}{n}),$
with ${x}'<0$,
$f(\tilde{\mbox\xx}_{1},\tilde{\mbox\vv}_{1})=({x}'',\frac{1}{n})$,
with ${x}''>0$. Since  $\cal A$ is connected and  $f$ is
continuous, by the intermediate value Theorem there exists
$({\mbox\xx}_{n},{\mbox\vv}_{n})\in\cal A$ such  that
$f({\mbox\xx}_{n},{\mbox\vv}_{n})=(0,\frac{1}{n})$.
Hence we have a sequence
$\{\,({\mbox\xx}_{n},{\mbox\vv}_{n})\,\}\in\cal A$ such  that
${\mbox\rv}(t_n,{\mbox\xx}_{n},{\mbox\vv}_{n})=(0,\frac{1}{n})$,
for some $t_n \in (0, T_{h_0})$ (in fact $t_n
=t_2({\mbox\xx}_{n},{\mbox\vv}_{n})$). Since $\cal A$ is compact
and $t_n \in [0, T_{h_0}],$  there exists  a subsequence $(t_m,
{\mbox\xx}_{m},{\mbox\vv}_{m})$ of
$(t_n,{\mbox\xx}_{n},{\mbox\vv}_{n})$ such  that  $(t_m,
{\mbox\xx}_{m},{\mbox\vv}_{m})\rightarrow (\tau,
\bar{\mbox\xx},\bar{\mbox\vv})\in [0, T_{h_0}]\times \cal A.$
Consequently $\lim_{m\rightarrow +\infty}{\mbox\rv}(t_m,
{\mbox\xx}_{m},{\mbox\vv}_{m})= {\mbox\rv}(\tau,
\bar{\mbox\xx},\bar{\mbox\vv})$. On the other hand,
${\mbox\rv}(t_m,
{\mbox\xx}_{m},{\mbox\vv}_{m})=(0,\frac{1}{m})\rightarrow
(0,0)$, which implies that ${\mbox\rv}(\tau,
\bar{\mbox\xx},\bar{\mbox\vv}, 1)=(0,0).$ Therefore  ${\mbox\rv}:
[0,\tau\,]\rightarrow \R^2,$ $\,\,{\mbox\rv}(t)={\mbox\rv}(t,
\bar{\mbox\xx},\bar{\mbox\vv})$ satisfies (1), (2), (3) of the
definition of an essential part of a symmetric
figure eight orbit given in the introduction.

\begin{figure}[!htb]
 \centering
 \includegraphics[scale=0.45]{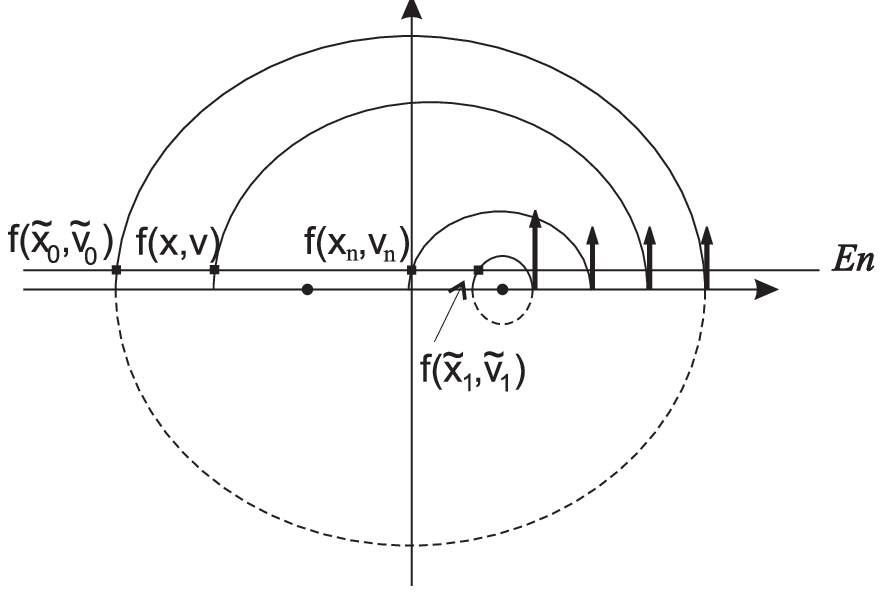}
 \hspace{2cm}
 \includegraphics[scale=0.8]{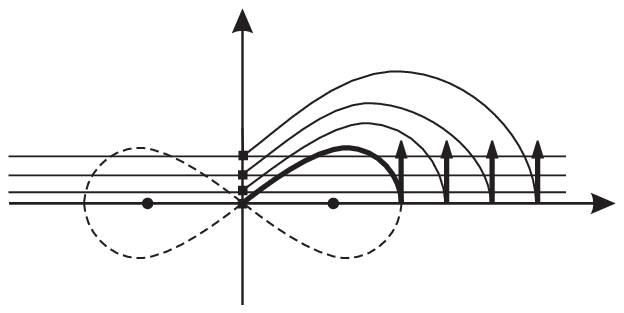}
 \caption{\scriptsize{Obtaining  the   solution  ${\mbox\rv}(t)$. }}
 \end{figure}

To finish the proof of the Theorem we show that
${\mbox\rv}(t)={\mbox\rv}( t, \bar{\mbox\xx},\bar{\mbox\vv})$
satisfies also property (4). Let $(x_m (t), z_m (t))={\mbox\rv}(t,
{\mbox\xx}_{m},{\mbox\vv}_{m}).$ Note that
$z_m (t)\geq 0,$ for all $t\in [0,t_m]$ (because $t_m =t_2
({\mbox\xx}_{m},{\mbox\vv}_{m})$). Then $\bar{z}(t)\geq 0,$ for
all $t\in [0,\tau\,],$ where $(\bar{x}(t),\bar{z}(t))={\mbox\rv}(
t, \bar{\mbox\xx},\bar{\mbox\vv})$. If
$\bar{z}(t')=0$ for some $t' \in(0,\tau)$, we have $\dot{\bar
z}(t')=0$ ($t'$ is a point of local minimum), that is,  $(\bar x,
\bar z)$ is tangent to the  $x$-axis at $t=t'$. It follows that
$(\bar{x}(t),\bar{z}(t))$ is  contained in the $x$-axis, for all
$t\in (0,\tau)$ (because the $x$-axis is an invariant subspace).
This is a contradiction. Then  $\bar{z}(t)>0$,
$\,t\in(0,\tau)$.
\CaixaPreta\\

\vspace{0,2cm}

\noindent{\bf Proof of Claim  \ref{1.8}.}
First we prove that we can choose ${\mbox{\sv}}_{0}(t)$
with small negative energy and small velocity.
Let ${\mbox{\rv}}_{0}(t)$ and ${\mbox{\rv}}_{\epsilon}(t)$ be as in Remark \ref{2.3}. Recall that
$\mbox{\sv}_{0}(t)=\frac{1}{\epsilon}\,{\mbox{\rv}}_{\epsilon}(\epsilon^{{3/2}}t)$.
Let $h, \bar h$ be  the energy of $\mbox{\rv}_{0}$ and ${\mbox{\rv}}_{\epsilon}$
respectively. We have $h<0$.
If ${\mbox{\vv}}_{\epsilon}$ is close to ${\mbox{\vv}}_{0}$ we have
that $\bar h$ is close to $h$. In particular, we can choose
${\mbox{\rv}}_{\epsilon}$ with
 energy  $\bar h= E({\mbox{\rv}}_{\epsilon}(t))<0$.
Note that $\tilde{\mbox\vv}_0 = {\dot{\mbox{\sv}}}_0 (0)=\epsilon\,
{\dot{\mbox{\rv}}}_{\epsilon} (0)=\epsilon \, {\mbox{\vv}}_\epsilon.$
Hence we can assume
$\| \tilde{\mbox\vv}_0\|=| {\tilde{v}}_0 | \leq \epsilon (2\|
{\mbox{\vv}}_0 \|).$
Also,
{\small $
E(\mbox{\sv}_{0}(t))=\frac{1}{2}\epsilon\,
\|\dot{{\mbox{\rv}}_{\epsilon}}(0)\|^{2}+V(\frac{1}{\epsilon}{\mbox{\xx}}_{0},1)
=\frac{1}{2}\epsilon\,\|{\mbox{\vv}}_\epsilon \|^{2}+\epsilon
V({\mbox{\xx}}_{0},\epsilon)
= \epsilon\,E({\mbox{\rv}}_{\epsilon}(t))\,=\,\epsilon\bar h <0,
$}
where $V({\mbox{\xx}},\epsilon)$ is as in section 2.1 (for the second
equality see Lemma 1.1).
This shows  that we can choose ${\mbox{\sv}}_0$ with
small velocity $\tilde{\mbox{\vv}}_0$ and small negative energy.

We now deal with ${\mbox{\sv}}_{{1}}(t)$. Let
${\mbox{\rv}}_\epsilon$, ${\mbox{\vv}}_\epsilon =(0, v_\epsilon)$,
${\mbox{\vv}}_0 =(0, v_0)$ as in Remark \ref{r2.8}. Recall that
${\mbox{\sv}}_{{1}}(t)=
\epsilon\,{\mbox{\rv}}_{{\epsilon}}(\frac{1}{\epsilon}t)$. Since
we can choose ${\mbox{\vv}}_\epsilon$ close to ${\mbox{\vv}}_0$,
we can assume that $0<a<\|{\mbox{\vv}}_\epsilon\|<b$, for some
constants $a, b$ independent of $\epsilon$. A simple calculation
shows that ${\dot{\mbox\sv}}_{{1}} (0)=
{\dot{\mbox{\rv}}}_{{\epsilon}} (0) ={{\mbox{\vv}}_\epsilon}.$ We
have $ E({\mbox\sv}_{1}(t))=\frac{1}{2}\,
\|{\dot{\mbox\sv}}_{{1}}(0)\|^{2}+
V({\mbox{\sv}}_{1}(0))<\frac{1}{2}b^2 +V({\mbox{\sv}}_{1}(0))$ and
that  $dist(\,{{\mbox\sv}}_1 (0),{\cal C})=\epsilon
\|{\mbox{\xx}}_0\|$, where ${\mbox{\xx}}_0$ (fixed) is as in Remark
\ref{r2.8}. Since $V({\mbox{\xx}})\rightarrow -\infty,$ when
$\,dist \,({\mbox{\xx}} ,{\cal C})\rightarrow 0$ (see \cite{AO1},
Lemma 2.4), we have that for any $n>0$ we can choose $\epsilon$
sufficiently  small  such  that $E(\mbox{\sv}_{1}(t))<-n$.  In
particular, we can choose $\mbox{\sv}_{0}(t)$ and
$\mbox{\sv}_{1}(t)$ such that $E(\mbox{\sv}_{1}(t))<
E(\mbox{\sv}_{0}(t))$. Also, since
$\|{\mbox\vv}_{1}\|=\|{\mbox\vv}_{\epsilon}\|>a$, and
$\|{\mbox\vv}_{0}\|$ can be chosen small, it
 follows that we can choose
$\|\tilde{\mbox\vv}_{0}\|<\|\tilde{\mbox\vv}_{1}\|.$ This  proves
the claim. \CaixaPreta

\vspace{0,4cm}

\noindent{\bf Proof of Claim  \ref{2.8}.} It follows  from claim  \ref{1.8} and from the
fact that $V(\tilde{\mbox{\xx}}_{1})\leq V(\mbox\xx)\leq
V(\tilde{\mbox{\xx}}_{0})$, for all
$\mbox\xx\in[\tilde{\mbox{\xx}}_{1},\tilde{\mbox{\xx}}_{0}\,].$
\CaixaPreta

\vspace{0,4cm}

\noindent{\bf Proof of Claim  \ref{4.8.8}.} Suppose that there is no such $\xi>0.$
Then there exists a sequence $(\mbox\xx_{n}, \mbox\vv_{n})$ such
that ${\tilde z}(\mbox\xx_{n},\mbox\vv_{n})\leq\frac{1}{n}$ and
$(\mbox\xx_{n},\mbox\vv_{n})$ converges to some
$(\bar{\mbox\xx},\bar{\mbox\vv})\in \cal A$.
Hence
${\mbox{\rv}}(t,\mbox\xx_{n},\mbox\vv_{n})\rightarrow
{\mbox{\rv}}(t,\bar{\mbox\xx},\bar{\mbox\vv})$ and
$\dot{\mbox{\rv}}(t,\mbox\xx_{n},\mbox\vv_{n})\rightarrow
\dot{\mbox{\rv}}(t,\bar{\mbox\xx},\bar{\mbox\vv}),$ which
implies that
$z_{n}(t)\rightarrow \bar{z}(t)$ and
$\dot{z}_{n}(t)\rightarrow \dot{\bar z}(t)$
 where $(x_{n}(t),
z_{n}(t))={\mbox{\rv}}(t,\mbox\xx_{n},\mbox\vv_{n})$ and
$(\bar{x}(t),\bar{z}(t))={\mbox{\rv}}
(t,\bar{\mbox\xx},\bar{\mbox\vv}).$

Since $\dot{\bar z} (t)>0$, for
$t\in [0, \tilde{t}(\bar{\mbox\xx},\bar{\mbox\vv}))$,
 we have that  there exists  $\chi>0$ such  that $\dot{\bar{z}}
(t)\geq \chi,$ for all $t\in
[0,\frac{1}{2}\tilde{t}(\bar{\mbox\xx},\bar{\mbox\vv})]$.
Hence for $n$ sufficiently  large
$\dot{z}_{n} (t) \geq\frac{\chi}{2}>0,\,\,\,\,\mbox{for
all}\,\,\,t\in
[0,\frac{1}{2}\tilde{t}(\bar{\mbox\xx},\bar{\mbox\vv})]$.
Also, since $z_{n}\rightarrow \bar{z}$ we have
 for $n$ sufficiently  large,
{\small $z_{n} (\frac{1}{2}\tilde{t}(\bar{\mbox\xx},\bar{\mbox\vv}))
\geq\frac{1}{2}
z(\frac{1}{2}\tilde{t}(\bar{\mbox\xx},\bar{\mbox\vv})) >0.$}
Note that ${\mbox{\rv}}(t,\mbox\xx_{n},\mbox\vv_{n})$   is
defined in
$[0,\frac{1}{2}\tilde{t}(\bar{\mbox\xx},\bar{\mbox\vv})]$, for
$n$ sufficiently  large.
By the definition of $\tilde t$, we have that
$\tilde{t}(\mbox\xx_{n},\mbox\vv_{n})>\frac{1}{2}\tilde{t}
(\bar{\mbox\xx},\bar{\mbox\vv})$ and, since $z_n(t)$ is increasing in
$[0,\tilde{t}({\mbox\xx},{\mbox\vv})]$ we get
{\small ${\tilde z}(\mbox\xx_{n},\mbox\vv_{n})=z_n(\tilde{t}(\mbox\xx_{n},\mbox\vv_{n}))
> {z}_{n}(\frac{1}{2}\tilde{t}(\bar{\mbox\xx},\bar{\mbox\vv}))
\geq \frac{1}{2}
{z}(\frac{1}{2}\tilde{t}(\bar{\mbox\xx},\bar{\mbox\vv}))=\eta >0
$} for all
$n$ sufficiently  large, a contradiction.
  \CaixaPreta

\vspace{0,5cm}
 Up to now have proved only the existence of one
symmetric figure eight orbit.
We show now how to obtain  infinitely many symmetric figure eight orbits.
First, note that the symmetric figure eight orbit constructed above has
its initial values $(\mbox\xx, \mbox\vv)$ in the set $\cal A$.
This set depends only on the choice of the initial values of the
solutions ${\mbox\sv}_0$ and ${\mbox\sv}_1$. These initial values
are $\tilde{\mbox\xx}_{0} =(\tilde{x}_0, 0),$
$\,\tilde{\mbox\vv}_{0}=(0, \tilde{v}_0)$ for ${\mbox\sv}_0$ and
$\tilde{\mbox\xx}_{1} =(\tilde{x}_1, 0),$
$\,\tilde{\mbox\vv}_{1}=(0, \tilde{v}_1)$ for ${\mbox\sv}_1$.
Its clear from the proof of claim \ref{1.8} that we can choose
$\tilde{v}_0 >0$ as small as we want and $\tilde{x}_1 >1$ as close
to $1$ as we want.
Hence are can choose a sequence $\{{\cal A}_n\}$ of disjoint sets,
all satisfying the statements of claims \ref{1.8} and \ref{2.8} (see
figure below).
Since for each ${\cal A}_n$ we have a symmetric figure eight orbit with
initial values on ${\cal A}_n$, we obtain in this way infinitely
many geometrically distinct figure eight orbits.

\begin{figure}[!htb]
\hspace{2cm}
 \begin{minipage}[b]{0.23\linewidth}
 \includegraphics[width=\textwidth]{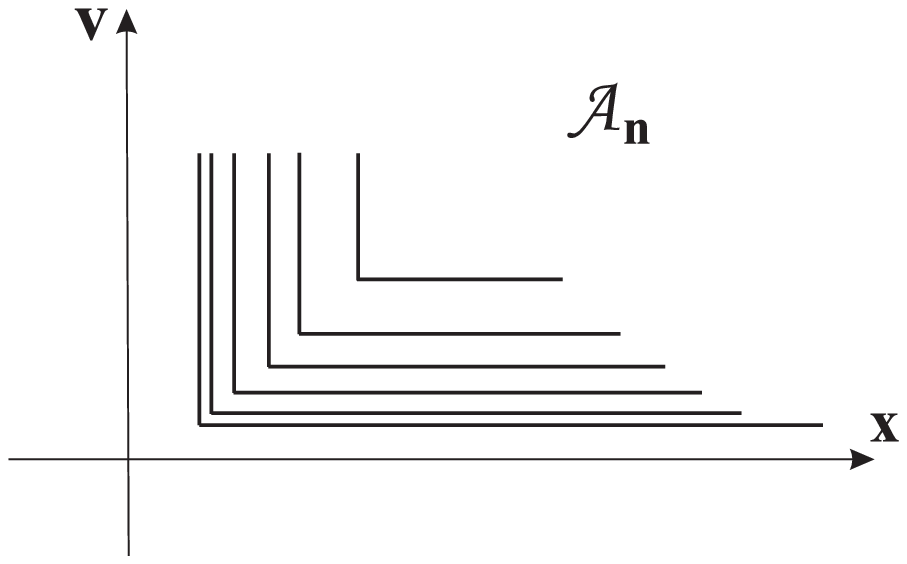}
\caption{\scriptsize{A sequence of disjoint sets ${\cal A}_n$.}}
 \end{minipage}  \hfill
\begin{minipage}[b]{0.23\linewidth}
 \includegraphics[width=\linewidth]{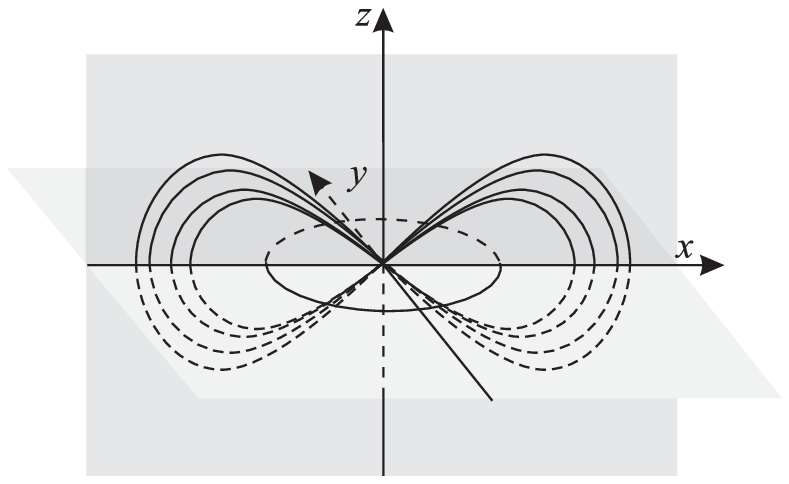}
 \caption{\scriptsize{Infinitely many figure eight orbits.}}
 \end{minipage} \hspace{2cm}
\end{figure}

\section{Proof of the Corollary to Theorem B: Symmetric figure eight orbits in the symmetric Euler problem.}

Consider the $xy$-plane and two
fixed centers with mass $M$ located at   $(\rho, 0)=\rho e_1$,
$\,(-\rho, 0)=-\rho e_1,$ $\rho>0.$ Consider a particle  $P$
moving in the $xy$-plane under the influence of the gravitational
attraction  induced by the two  fixed centers.
The potential of this  problem is given by

{\footnotesize $$U(\mbox{\rv})=-\frac{M}{\|
\mbox{\rv}+\rho e_1\|} - \frac{M}{\| \mbox{\rv}-\rho e_1\|}.$$}
Write $U( \mbox{\rv},\rho,M)=-\frac{M}{\| \mbox{\rv}+\rho
e_1\|} - \frac{M}{\| \mbox{\rv}-\rho e_1\|}$ to express the fact that
$U$ depends on the mass $M$ and the distance $\rho$.
$\,U$ satisfies the same properties as $V$ (section 1).
For  $\rho $ small $\,U( \mbox{\rv},\rho,M)$ is a symmetric
perturbation of the Kepler potential
$U(\mbox{\rv},0,M)=-\frac{2M}{\|\mbox{\rv}\|}$.
We can then
obtain (using the results in \cite{AO}), as in the case of the fixed homogeneous circle problem, periodic
symmetric orbits $\mbox{\sv}_0 (t)$ far from the two fixed
centers. Moreover, a similar argument as the one used in the proof
of claim  \ref{1.8}, shows  that we can choose  $\mbox{\sv}_0 (t)$
with small negative energy, and  also such  that $\|\dot{\mbox\sv}
_0 (t)\|$ is small.
\begin{wrapfigure}[11]{r}{3.8cm}
 \centering
 \includegraphics[width=3.5cm]{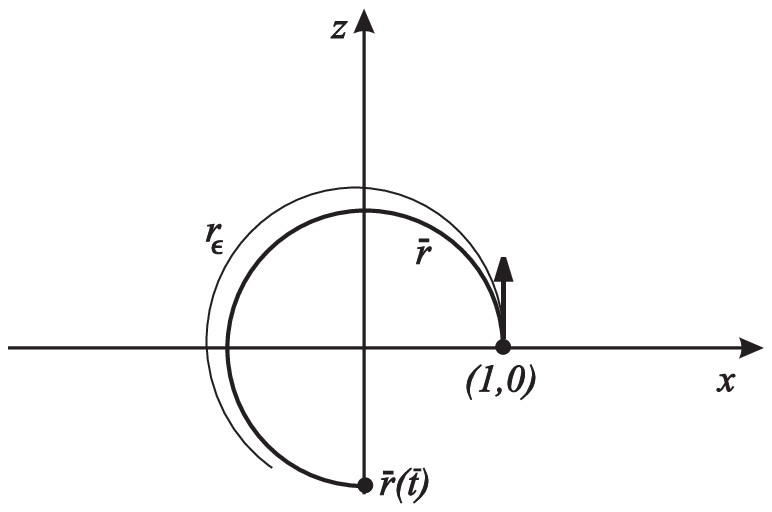}
 \caption{\scriptsize{Solutions $\bar{\mbox\rv}$ and ${\mbox\rv}_\epsilon$. }}
 \end{wrapfigure}
Let
{\footnotesize $${\cal U}({\mbox\rv},\epsilon)=
-\frac{M}{\|{\mbox\rv}\|}-\frac{M}{\|{\mbox\rv}+2\frac{1}{\epsilon}
e_1\|}= -\frac{M}{\|{\mbox\rv}\|}-\frac{\epsilon
M}{\|\epsilon{\mbox\rv}+2e_1\,\|}.$$}

\noindent Note that $U({\mbox\rv},\frac{1}{\epsilon},M)$ is obtained by a translation
from ${\cal U}({\mbox\rv},\epsilon)$.
 Note
also that for $\epsilon$ small we can consider ${\cal
U}({\mbox\rv},\epsilon)$ as a perturbation of Kepler potential
${\cal U}({\mbox\rv},0)=-\frac{M}{\|{\mbox\rv}\|}.$ Let
$\bar{\mbox\vv}=(0,\bar v\,),\,\,\bar v >0$, be a velocity  such
that $\bar{\mbox\rv}(t)$ is a circular solution  of the   Kepler
problem $\stackrel{..}{\mbox\rv}\,\,=-\nabla {\cal
U}({\mbox\rv},0),$ with $\bar{\mbox\rv}(0)=(1,0)=e_1$ and
$\dot{\bar{\mbox\rv}}(0)=\bar{\mbox\vv}$. Note that
$\bar{v}=\sqrt{M}.$
Let  also  $\bar t >0$ be such  that $\bar{\mbox\rv}(\bar
t)=(0,-1)$, $\bar t>0$ is the first (positive) time such that
$\bar{\mbox\rv}({\bar{t}})=(0,-1)$. Note that $\bar{\mbox\rv}(t)$,
$\,\,0\leq t\leq \bar t$, intersects transversally the negative $x$-axis
in a single point.
Let ${\mbox\rv}_\epsilon (t)= ({x}_\epsilon
(t),{y}_\epsilon (t))$ be a solution of
$\stackrel{..}{\mbox\rv}\,\,=\nabla {\cal
U}({\mbox\rv},\epsilon),$
with ${\mbox\rv}_\epsilon (0)=e_1,\,\dot{\mbox\rv}_\epsilon (0)=\bar {\mbox\vv}.$
By Proposition  1.4 of \cite{AO}, for $\epsilon$ sufficiently small,
 ${\mbox\rv}_\epsilon (t),\,\,0\leq t\leq \bar t,$
intersects   transversally the negative $x$-axis in exactly one
point  ${\mbox\rv}_\epsilon (t_\epsilon)$ and $y_\epsilon
(t)>0,\,\,0<t<t_\epsilon.$
Since ${\mbox\rv}_{\epsilon}$ is close
to $\bar{\mbox{\rv}}$ we can
assume that $-2< x_\epsilon (t_\epsilon).$
Let $\tilde{\mbox\rv} _\epsilon (t):= {\mbox\rv}_\epsilon (t)
+\frac{1}{\epsilon}e_1.$ Then $\tilde{\mbox\rv}_\epsilon (t)$  is
a solution of $\stackrel{..}{\mbox\rv}\,\,=\nabla
U({\mbox\rv},\frac{1}{\epsilon},M).$ Let ${\mbox\sv}^\epsilon (t)
:={\epsilon}\,\tilde{\mbox\rv} _\epsilon (\epsilon^{-3/2} t).$
It follows from   Corollary \ref{1.3.2} that ${\mbox\sv}^\epsilon (t)
=(x^\epsilon (t), y^\epsilon (t))$ is a solution of
$\stackrel{..}{\mbox\rv}\,\,=\nabla U({\mbox\rv},1,M),$
 and it is easy to verify  that ${\mbox\sv}^\epsilon
(0)=(1+{\epsilon},0),$ $\,{\mbox\sv}^\epsilon
(\epsilon^{3/2}t_\epsilon)=(x^\epsilon ,0),\,0<x^\epsilon <1,$
$\,y^\epsilon (t)>0,\,\,0<t<\epsilon^{3/2}t_\epsilon,$ and
$\dot{\mbox\sv}^\epsilon (0)=\frac{1}{\sqrt \epsilon}
\bar{\mbox\vv}.$ Then

{\footnotesize $$
E({\mbox\sv}^\epsilon (t))\,=\,E({\mbox\sv}^\epsilon (0)) \, =\,
\frac{1}{2\epsilon} \,\|\bar{\mbox\vv}\|^2 + U({\mbox\sv}^\epsilon
(0),1,M)\, =\, \frac{M}{2\epsilon}
-\frac{M}{\epsilon}-\frac{M}{2+{\epsilon}}\,
=\,  -\frac{M}{2\epsilon}\Bigl(1+\frac{2\epsilon}{2+\epsilon}\Bigr).$$}

Hence for $\epsilon$ sufficiently  small we can choose
${\mbox\sv}_1 (t)={\mbox\sv}^\epsilon (t)$ with large and negative energy.
Note  also that
$\|\dot{\mbox\sv}_1
(0)\|=\frac{1}{\sqrt\epsilon}\,\|\bar{\mbox\vv}\|= \frac{\sqrt
M}{\sqrt\epsilon}$ is large.
The solutions  ${\mbox\sv}_0 (t)$ and  ${\mbox\sv}_1 (t)$ constructed above
satisfy the same properties as the solutions  ${\mbox\sv}_0 (t)$ and
${\mbox\sv}_1 (t)$ constructed at the beginning of the proof of Theorem B.

For the case of the symmetric Euler problem, Claim
\ref{1.8} follows from the choice of ${\mbox\sv}_0$ and ${\mbox\sv}_1$
above; the proofs of Claims  \ref{2.8} and \ref{4.8.8} are
identical. The proof of Proposition \ref{3.8} is similar, just estimate
$\frac{\partial}{\partial y} U(x,y) = y
\Bigl(\frac{M}{\|\mbox{\rv} -e_1 \|^3} +\frac{M}{\|\mbox{\rv} +e_1
\|^3}\Bigr)$, where $U(x,y)=U((x,y),1,M).$ But a simple calculation
shows  that if $\mbox\rv(t) =(x(t),y(t))$ has energy
less or equal $\delta<0,$ then

{\footnotesize $$\Bigl(\frac{M}{\|\mbox{\rv} -e_1 \|^3} +\frac{M}{\|\mbox{\rv} +e_1
\|^3}\Bigr)\geq \frac{2M}{(1+R_\delta )^3},$$}

\noindent where $R_\delta <+\infty$ is the radius of  Hill's region
with energy  $\delta<0$, that is,  $R_\delta= sup \,\{\,\|\mbox\rv
\|\,;\, U(\mbox\rv )\leq \delta \,\}$ (it can be verified directly that
$R_\delta \leq \frac{2}{-\delta} +1$).
Then    $y(t)$ satisfies
$|\stackrel{..}y\,|\geq A\,|y|,$
with $A=\frac{2M}{(1+R_\delta )^3}.$
The rest of
the proofs of Proposition \ref{3.8} and of  Theorem  B for
this case are similar. \vspace{0,2cm}

As in the case of the fixed homogeneous circle, we can prove that
for the symmetric Euler problem we also have infinitely many
symmetric figure eight orbits. \CaixaPreta

\vspace{0.2cm}

\section{Proof of Theorem C: Periodic Spiral Orbits.}

In this section we fix the density $\lambda$ of the circle and denote by
$V(\mbox{\rv},\frac{1}{\epsilon})=V(x,y,z; \frac{1}{\epsilon})$ the potential at the point
 $\mbox{\rv}=(x,y,z)\in\R^3$ induced by
the fixed homogeneous circle with density $\lambda$, contained in the
$xy$-plane, centered at the origin
and with radius $\frac{1}{\epsilon}$. Since $V$ is invariant by rotations around the
$z$-axis we can reduce our problem in a canonical way to a problem with
two degrees of freedom. The Lagrangian in
cylindrical coordinates can be written as:
$L(r,\varphi,z,\dot r,\dot\varphi,\dot z)=\frac{1}{2}({\dot
r}^{2}+r^{2}\dot\varphi ^{2}+\dot z ^{2})-V_{\epsilon}(r,z)$,
where
$V_{\epsilon}(r,z)=V(rcos\varphi , r sin\varphi,z,\frac{1}{\epsilon})$.
It is then straightforward to verify that the system
$\stackrel{..}{\mbox{\rv}}\,=-\nabla V({\mbox{\rv}};\frac{1}{\epsilon})$
in these coordinates is given by:

{\footnotesize \begin{equation}
\left\{ \begin{array}{l}
\stackrel{..}{r}\,  =\displaystyle\frac{\frac{K^{2}}{\epsilon^{2}}}{
r^{3}}
-\frac{\partial V_{\epsilon}}{\partial r}
\left( {r},{z}\right)\\  \\
\stackrel{..}{z}\,  =-\displaystyle\frac{\partial V_{\epsilon}}{\partial z}
\left( r,z\right)\\ \\
\dot\varphi(t)=\displaystyle\frac{K/\epsilon}{r^{2}}
\end{array}
\right.
\label{R4}
\end{equation}}
\noindent where $\frac{K}{\epsilon}$ is the (constant) angular momentum.

\vspace{0,3cm}

\noindent {\bf Remarks.}

\noindent (1) $K=0$ implies that $\varphi$ is constant.  Then the
particle  moves on the vertical plane  determined by the $z$-axis and by
the vector ($cos\varphi,sin\varphi,0$).

\noindent (2)  The first two equations of system (\ref{R4}) can be
rewritten  as $\stackrel{..}{{{\mbox{\rv}}}}\,  =-\nabla {\overline
V}\Bigl({\mbox{\rv}},\frac{1}{\epsilon}\Bigr),$ with
${\mbox{\rv}}=(r,z)$ and
{\footnotesize $$\overline V  \Bigl(r,z,\frac{1}{\epsilon}\Bigr)=
\frac{
\frac{K^{2}}{\epsilon^{2}}}{2 r^{2}}
+\,V \left( {r},0,{z},\frac{1}{\epsilon}\right).$$}

\vspace{0,2cm}

Note that, if $(r(t),z(t))$ is a solution of the first two equations of (\ref{R4}), defining
$\varphi(t)=\displaystyle\int_{0}^{t}\frac{K/\epsilon \,ds}{r^{2}(s)},$ we have
that $(r(t),\varphi(t),z(t))$ is a solution of (\ref{R4}).
Then $(r(t)\,cos\varphi(t),r(t) \, sin\varphi(t),z(t))$
is a solution of $\stackrel{..}{\mbox{\rv}} =-\nabla V({\mbox{\rv}};\frac{1}{\epsilon})$.

\vspace{0,2cm}

{\lem Let $(r(t),z(t))$ be a periodic solution  of the first two equations of
(\ref{R4}) with  period
$\tau$ and {\footnotesize$\varphi(t)=\displaystyle\int_{0}^{t}\frac{K/\epsilon \,ds}{r^{2}(s)}$.}
Then ${\mbox{\rv}}(t)=(r(t)\,cos\varphi(t),r(t) \, sin\varphi(t),z(t))$
is a periodic solution  of
$\stackrel{..}{\mbox{\rv}}\,=-\nabla V({\mbox{\rv}};\frac{1}{\epsilon})$ if and only if
{\footnotesize$\,\displaystyle\frac{\varphi(\tau)}{2\pi}$} is rational.
\label{4.4.2}}

\vspace{0,25cm}

\noindent {\bf Proof.} Suppose
{\small$\frac{\varphi(\tau)}{2\pi}=\frac{p}{q}$},
 $p,q\in \Z$. Let $R_\theta$ denote the rotation about
the $z$-axis by an angle $\theta$. Note that $(R_\theta)^l= R_{l\theta}$
for every $l\in\Z$.
 Since $(r(t),z(t))$ is
$\tau$-periodic and {\small$\dot{\varphi}=\frac{K/\epsilon}{r^2}$} a
direct calculation shows that:
${\mbox{\rv}}(\tau)\,=\, R_{\varphi(\tau)} {\mbox{\rv}}(0)$ and
 $\dot{{\mbox{\rv}}}(\tau)\,=\, R_{\varphi(\tau)} \dot{{\mbox{\rv}}}(0).$
This implies
${\mbox{\rv}}(q\tau)\,=\, R_{q\varphi(\tau)} {\mbox{\rv}}(0)$ and
 $\dot{{\mbox{\rv}}}(q\tau)\,=\, R_{q\varphi(\tau)} \dot{{\mbox{\rv}}}(0).$
Since $q\varphi(\tau)= 2\pi p,$  then ${\mbox{\rv}}(q\tau)\,=\,
{\mbox{\rv}}(0),\,\, \dot{{\mbox{\rv}}}(q\tau)\,=\,\dot{{\mbox{\rv}}}(0).$
Hence ${\mbox{\rv}}$ is $q\tau$ periodic.

Conversely, suppose  that    $(r(t)\,cos\varphi(t),r(t) \,
sin\varphi(t),z(t))$ is $\tau_{0}$-periodic, for some $\tau_{0}$.
Then, $z(t)$ and $r(t)=\sqrt{
(r(t)^{2}\,cos^{2}\varphi(t)+r(t)^{2} \, sin^{2}\varphi(t)}$ are
$\tau_{0}$-periodic.
Let $\tau$ be the minimal period of
$(r(t),z(t))$.
From the definition of $\varphi$ it follows that
$\varphi(k\tau)=k\varphi(\tau)$, for all $k\in\Z$.
It also follows  that $\tau_{0}=n\tau$.
Hence, since  $r(t)>0$, for all $t$, and $r(t)$  is
$\tau_{0}$-periodic, we have
$cos\varphi(t+n\tau)=cos\varphi(t+\tau_{0})=cos\varphi(t)$. In the
same way $sin\varphi(t+n\tau)=sin\varphi(t+\tau_{0})=sin\varphi(t)$,
which implies that $\varphi(t+n\tau)-\varphi(t)=2\pi p$, for some
$p\in\Z$.
Evaluating in $t=0$, we have $n\varphi(\tau)=\varphi(n\tau)=2\pi p$,
and therefore  {\footnotesize$\displaystyle\frac{\varphi(\tau)}{2\pi}=\frac{p}{n}$}.
\CaixaPreta\\

\vspace{0,2cm}

Recall that in section 2.2 we denoted by $W(x,z;\epsilon)$ the potential at the point $(x,z)$
(in the $xz$-plane) induced by a circle in the $xy$-plane of radius $\frac{1}{\epsilon}$ centered
at $(\frac{1}{\epsilon},0,0)$. Now, for $K\in\R$  consider the system:
{\footnotesize \begin{equation}
\left\{ \begin{array}{l}
\stackrel{..}{r}\,  =\displaystyle\frac{
\frac{K^{2}}{\epsilon^{2}}}{\left(
r-\frac{1}{\epsilon}\right)^{3}}
-\frac{\partial W}{\partial r}
\left({r},{z};{\epsilon}\right)\\  \\
\stackrel{..}{z}\,  =-\displaystyle\frac{\partial W}{\partial z}
\left( r,z;{\epsilon}\right)
\end{array}
\right.
\label{R2}
\end{equation}}

\noindent{\bf Remarks.}

\noindent (1) System (\ref{R2}) can be
rewritten as  $\stackrel{..}{{{\mbox{\rv}}}}\,  =-\nabla {\overline
W}\Bigl({\mbox{\rv}};{\epsilon}\Bigr),$ where
{\footnotesize
$$\overline W \Bigl({\mbox{\rv}};{\epsilon}\Bigr)=
\overline W \Bigl(r,z;{\epsilon}\Bigr)=  \frac{
\frac{K^{2}}{\epsilon^{2}}}{2\left(
r-\frac{1}{\epsilon}\right)^{2}}
+\,W \left({r},{z};{\epsilon}\right).$$}

\noindent (2) Note that $\overline V$ is obtained from $\overline
W$ by a translation. Indeed $\overline V
\Bigl(x-\frac{1}{\epsilon},z,\frac{1}{\epsilon}\Bigr)=\overline
W \Bigl(x,z;{\epsilon}\Bigr).$

\vspace{0,25cm}

Recall that a key ingredient in the proof of Theorem A, part (ii), was the fact that for $K=0$
system \ref{R2} is a perturbation of the
infinite homogeneous straight wire problem. The next Lemma says that the same is true for any $K$.

{\lem For every $K$ the system given by (\ref{R2})
behaves as a perturbation of the
infinite homogeneous straight wire problem.}

\noindent{\bf Proof.} The Lemma follows from Proposition \ref{2a} and the fact that  $lim_{\epsilon\rightarrow 0}\frac{
\frac{K^{2}}{\epsilon^{2}}}{\left(
x-\frac{1}{\epsilon}\right)^{3}}=0$ uniformly on compacts. \CaixaPreta \\

\vspace{0,2cm}

We will use the following notation. Let $U\subset\{(x,0,z); \,x>0\}\subset\R^{2}$.
 Define the rotation of $U$
about the $z$-axis: $rot\,U= \{\,(x\, cos\varphi,x
\,sin\varphi,z);\,\,\varphi\in\R,\,\,\,(x,z)\in U\,\}.$ In
cylindrical  coordinates, we have $rot\,U= \{(r,\varphi,z);\,
(r,z)\in U\}.$ The next result, together with some rescaling, will
imply Theorem C.

{\prop
Let $C$ be a circle  in the $xz$-plane,
with center at the origin and
let $\,U$ be an open bounded set that contains $C$, of the form
$C\subset U\subset (\R^{2}-\{ (0,0)\})$.
Let $K\neq 0.$
Then  for each
$\epsilon_{0}>0$ there exist  $\epsilon$,
$0<\epsilon<\epsilon_{0}$, and a periodic solution
$\mbox{\rv}_{\epsilon}(t)=(r_{\epsilon}(t)
cos\varphi_{\epsilon}(t), r_{\epsilon}(t)sin
\varphi_{\epsilon}(t),z_{\epsilon}(t))$
of
$\,\stackrel{..}{{{\mbox{\rv}}}}\,  =-\nabla {
V}\Bigl({\mbox{\rv}},\frac{1}{\epsilon}\Bigr)$ in $\,\,rot
\left(U+\{(\frac{1}{\epsilon},0)\}\right)$ with angular momentum
$K/\epsilon$. Moreover, the trace of $(r_{\epsilon}(t),z_{\epsilon}(t))$
is a simple closed curve, symmetric  with respect to the $x$-axis,
and encloses the fixed homogeneous circle. \label{4.4.3}}

\vspace{0,24cm}

\noindent {\bf Proof.} Let $C$ be a circle in the $rz$-plane
with center at the  origin and
 $\,U$ be an open bounded set that contains $C$ of the form $C\subset U\subset
 \R^{2}-\{(0,0)\}.$ Without loss of generality we can assume $K>0$.
Let $\mbox{\xx}_{0}$  be the unique point in $C\cap(\mbox{positive}\,r-\mbox{axis}).$
We will denote by ${\bar{\mbox\rv}}_{{}_{\mbox\vv},{}_{\mbox{$\epsilon$}}}(t)=
({\bar{r}}_{{}_{\mbox\vv},{}_{\mbox{$\epsilon$}}}(t),{\bar{z}}_{{}_{\mbox\vv},{}_{\mbox{$\epsilon$}}}(t))$
a solution of
 $\,\stackrel{..}{{{\mbox{\rv}}}}\,  =-\nabla {\overline W}({\mbox{\rv}};{\epsilon})$,
with ${\bar{\mbox\rv}}(0)=\mbox{\xx}_{0}$ and
$\dot{\bar{\mbox\rv}}(0)=\mbox{\vv}.$

Let ${\bar{\mbox\rv}}_{{}_{\mbox\vv_{0}},{}_{\mbox{$0$}}}(t)$ be a circular solution
of $\,\stackrel{..}{{{\mbox{\rv}}}}\,  =-\nabla {\overline W}({\mbox{\rv}};0)$ whose trace is  $C$.
Consider     $[0,\bar t\,]$, a time interval in which
${\bar{\mbox\rv}}_{{}_{\mbox\vv _{0}},{}_{\mbox{$0$}}}(t)$ intersects
transversally the closed segment $E=\{\,(r,0)\,;\,-\infty<r\leq
0\,\}$ in a single point. By the continuous dependence of the solutions
(e.g. see Proposition 1.5 of \cite{AO})
we have that there is a $\delta>0$ such  that, if
$|\!|\mbox{\vv}-\mbox{\vv}_{0}|\!|<\delta,\,|\epsilon|<\delta,$
then
${\bar{\mbox\rv}}_{{}_{\mbox\vv},{}_{\mbox{$\epsilon$}}}(t)=({\bar
r}_{{}_{\mbox\vv},{}_{\mbox{$\epsilon$}}}(t),{\bar
z}_{{}_{\mbox\vv},{}_{\mbox{$\epsilon$}}}(t))$, $t\in [0,\bar
t\,]$, intersects   $E$ in a single point, and
$t(\mbox{\vv},\epsilon)$
is continuous, where $t(\mbox{\vv},\epsilon)$ is such that
${\bar{\mbox\rv}}_{{}_{\mbox\vv},{}_{\mbox{$\epsilon$}}}(t(\mbox{\vv},\epsilon))\in E$.
Define  $\tau(\mbox{\vv},\epsilon):=2t(\mbox{\vv},\epsilon)$ and
 $V_{\delta}=\{\,(1+s)\, \mbox{\vv}_{0};
\,|s|<\frac{\delta}{|\!|\mbox{\vv}_{0}|\!|}\,\}.$
We can assume

\noindent (1) $\tau(\mbox{\vv},\epsilon)$ is bounded,

\noindent (2) $|{\bar r}_{{}_{\mbox\vv},{}_{\mbox{$\epsilon$}}}(t)|$ is bounded,

\noindent (3) $0<\gamma <\tau(\mbox{\vv},\epsilon)$ for some constant $\gamma.$

Let ${\mbox{\rv}}_{{}_{\mbox\vv},{}_{\mbox{$\epsilon$}}}(t)=(
r_{{}_{\mbox\vv},{}_{\mbox{$\epsilon$}}}(t),
 z_{{}_{\mbox\vv},{}_{\mbox{$\epsilon$}}}(t))=
(\bar r_{{}_{\mbox\vv},{}_{\mbox{$\epsilon$}}}(t)+\frac{1}{{\epsilon}},\bar z
_{{}_{\mbox\vv},{}_{\mbox{$\epsilon$}}}(t))$.
Note that  ${\mbox{\rv}}_{{}_{\mbox\vv},{}_{\mbox{$\epsilon$}}}(t)$ is a
solution of $\,\stackrel{..}{{{\mbox{\rv}}}}\,  =-\nabla {\overline V}({\mbox{\rv}},\frac{1}{\epsilon})$.
 We can also choose $\delta$ small such that ${r}_{{}_{\mbox\vv},{}_{\mbox{$\epsilon$}}}>0,$
for all $\epsilon\in(-\delta,\delta)$.
Hence we have that ($ r_{{}_{\mbox\vv},{}_{\mbox{$\epsilon$}}},\varphi,
z_{{}_{\mbox\vv},{}_{\mbox{$\epsilon$}}}$) is solution of
(\ref{R4}), where
{\footnotesize $\varphi(t) :=
\varphi_{{}_{\mbox\vv},{}_{\mbox{$\epsilon$}}}(t)=\displaystyle\int_{0}^{t}
\frac{K/\epsilon \,\, ds}{{r}^{2}
_{{}_{\mbox\vv},{}_{\mbox{$\epsilon$}}}(s)}$}.
Define  $\Theta:V_{\delta}\times [0,\delta]\rightarrow\R$
by
{\footnotesize $$\Theta({\mbox\vv},\epsilon)=\left\{\begin{array}{l}
\displaystyle\frac{1}{2\pi}\int_{0}^{\tau({\mbox\vv},\epsilon)}\frac{K/\epsilon\, \,ds}
{ r ^{2}_{{}_{\mbox\vv},{}_{\mbox{$\epsilon$}}}(s)} = \frac{1}{2\pi}
\displaystyle\int_{0}^{\tau({\mbox\vv},\epsilon)}\frac{K/\epsilon \,\,ds}
{\left( \bar r
_{{}_{\mbox\vv},{}_{\mbox{$\epsilon$}}}(s)+\frac{1}{\epsilon}\right)^{2}},
\,\,\,\,\, \mbox{for }\,\,\epsilon\neq 0\\ \\
0,\,\,\,\,\,\,\,\mbox{for}\,\,\,\,\epsilon=0.
\end{array}
\right.$$}

From (1) and (2) above and from the fact that $\tau$ is continuous
it follows that
$\Theta $ is continuous in $V_{\delta}\times [0,\delta].$
Note that
$\Theta(\mbox{\vv},\epsilon)=\frac{\varphi_{{}_{\mbox\vv},{}_{\mbox{$\epsilon$}}}
(\tau(\mbox{\vv},\epsilon))}{2\pi}.$ Since $\tau(\mbox{\vv},\epsilon)\geq \gamma >0$
 (see (3) above) and we are assuming
$K>0$, we have that  $\,\Theta(\mbox{\vv},\epsilon)>0$ for
$\epsilon>0$. Also $\Theta(\mbox{\vv},0)=0,$ for all $\mbox{\vv}$.

Let $\epsilon_{0}>0$. By the Addendum to Theorem 0.2 of
\cite{AO} (taking $\epsilon_{0}>0$ even smaller,
if necessary) there exists a compact connected  set
${\cal V}\subset V_{\delta}\times [0,\epsilon_{0}] $ with
${\cal V}_{\epsilon}={\cal V}\cap (V_{\delta}\times\{\epsilon\})\neq\emptyset,$
for all
$\epsilon\in [0,\epsilon_{0}]$
such  that for
$(\mbox{\vv},\epsilon)\in {\cal V},$
$\,
\bar{\mbox{\rv}}_{{}_{\mbox\vv},{}_{\mbox{$\epsilon$}}}(t)$
is  a periodic   solution
of
$\,\stackrel{..}{{{\mbox{\rv}}}}\,  =-\nabla {\overline W}({\mbox{\rv}};{\epsilon})$
whose trace is a  simple closed curve
symmetric  with respect to the $r$-axis, and encloses
 the origin.

 Since ${\cal V}_{\epsilon_{0}}\neq\emptyset$
and ${\cal V}_{0}\neq\emptyset$, there exist
$\,(\mbox{\vv}_{\epsilon_{0}},\epsilon_{0}), (\mbox{\vv}_{0},0)\in
{\cal V}.$ Moreover,
$\Theta(\mbox{\vv}_{\epsilon_{0}},\epsilon_{0})>0$ and
$\Theta(\mbox{\vv}_{0},0)=0.$ Since $\cal V$ is connected we have that
$[0,\Theta(\mbox{\vv}_{\epsilon_{0}},\epsilon_{0})]\subset
\Theta({\cal V}).$ Then there exists a rational number
$\frac{p}{q}\in
(0,\Theta(\mbox{\vv}_{\epsilon_{0}},\epsilon_{0})]$ and hence
there exists  $(\mbox{\vv}_{\epsilon},\epsilon)\in \cal V$,
$\,0<\epsilon\leq\epsilon_{0}$, such  that
$\Theta(\mbox{\vv}_{\epsilon},\epsilon)=\frac{p}{q}$. From Lemma
\ref{4.4.2} it follows that
$(r_{{}_{\mbox\vv},{}_{\mbox{$\epsilon$}}}(t)
\,cos\varphi_{{}_{\mbox\vv},{}_{\mbox{$\epsilon$}}}(t),
r_{{}_{\mbox\vv},{}_{\mbox{$\epsilon$}}}(t)\,sin
\varphi_{{}_{\mbox\vv},{}_{\mbox{$\epsilon$}}}(t),
z_{{}_{\mbox\vv},{}_{\mbox{$\epsilon$}}}(t))$
is a periodic solution of
$\,\stackrel{..}{{{\mbox{\rv}}}}\,  =-\nabla {V}({\mbox{\rv}},\frac{1}{\epsilon})$,
 that satisfies the properties of the statement of the  Proposition \ref{4.4.3}.
 \CaixaPreta

\vspace{0,5cm}

\noindent{\bf Proof of Theorem C.}
In Proposition \ref{4.4.3} take $C$ with
radius $\frac{1}{2}$ and  $\,U=\{\,
\mbox{\pp}\in\R^2\,;\,\frac{1}{3}<\mbox\pp <1\,\}$.
Then  there exist $\epsilon$,
$0<\epsilon<\epsilon_{0},$ and a periodic solution
{\footnotesize $\mbox{\rv}_{\epsilon}(t)=(r_{\epsilon}(t)
\,cos\varphi_{\epsilon}(t), r_{\epsilon}(t)\,sin
\varphi_{\epsilon}(t),z_{\epsilon}(t))$} of
 $\,\stackrel{..}{{{\mbox{\rv}}}}\,  =-\nabla {V}({\mbox{\rv}},\frac{1}{\epsilon})$, with
angular momentum $K/\epsilon$. Moreover, the trace of
$(r_{\epsilon}(t),z_{\epsilon}(t))$ is a simple closed curve
symmetric with respect to the $r$-axis, and encloses the fixed
homogeneous circle. By the choice of $U$ we have $\frac{1}{3}<\,
dist(\mbox{\rv}_{\epsilon},{\cal C}_\epsilon )<1$.

Let {\small $r(t)=\epsilon
\,r_{\epsilon}\left(\frac{1}{\epsilon}t\right),$
$\,z(t)=\epsilon\, z_{\epsilon}\left(\frac{1}{\epsilon}t\right),$
$\,\varphi(t)=\varphi_\epsilon \left(\frac{1}{\epsilon}t\right),$}
and ${\mbox{\rv}}(t)=(r(t)\,cos\varphi(t),r(t) \,
 sin\varphi(t),z(t))$,
 i.e.
 $\mbox{\rv}(t)=\epsilon \,
\mbox{\rv}_{\epsilon}\left(\frac{1}{\epsilon}t\right)$. By
Corollary \ref{1.3.2} we have that $\mbox\rv (t)$  is a solution of
our original problem
$\,\stackrel{..}{{{\mbox{\rv}}}}\,=-\nabla {V}({\mbox{\rv}},1)\, =-\nabla {V}({\mbox{\rv}})$.
Also, by the  properties of $\mbox{\rv}_{\epsilon}$ we have
that $\mbox\rv (t)$ satisfies:
(1) \,$\frac{\epsilon}{3}< dist(\mbox\rv (t),{\cal C})
<\epsilon$,
(2) writing $\mbox\rv_\epsilon
(t)=(x_{\epsilon}(t),y_{\epsilon}(t),z_{\epsilon}(t))$ we see
that
{\small $\mbox{\rv}(t)=
 (\epsilon
\,x_{\epsilon}\left(\frac{1}{\epsilon}t\right), \epsilon\,
y_{\epsilon}\left(\frac{1}{\epsilon}t\right), \epsilon\,
z_{\epsilon}\left(\frac{1}{\epsilon}t\right))$} has angular momentum
$\epsilon\,({\dot y}_{\epsilon}x_{\epsilon}-{\dot
x}_{\epsilon}y_{\epsilon})=\epsilon \frac{K}{\epsilon}=K.$ \CaixaPreta


\section{Some Generalizations.}

In this section we indicate how the  existence of some periodic orbits remains
true if we add some other body $\cal B$ to the fixed homogeneous circle $\cal C$.
For a Lebesgue measurable ${\cal B}\subset \R^3$ and $\rho\in\R$, ${\cal B}_{\rho}$
denotes the set $\{\,\rho b\,;\,b\in{\cal B}\,\}$.
Let $\lambda$ be a positive finite measure on $\cal B$.
The (total) mass of $\cal B$ is
$M_{\cal B}=\int_{\cal B} \lambda =\lambda({\cal B}) > 0$.
For $\rho >0$, $\lambda_\rho$ denotes the measure on ${\cal B}_{\rho}$ induced
by $\lambda$, i.e. $\lambda_\rho (A_\rho)=\lambda(A)$, for all Lebesgue measurable $A\subset\cal B$.
For $\mbox{\rv}\notin \overline{\cal B}$, define  $V({\mbox{\rv}}, \rho, M)=\frac{M}{M_{\cal
B}}\int_{{\cal B}_{\rho}}
\frac{\lambda_{\rho}(u)}{\|{\mbox{\rv}}-u\|}\,du$.
It is
straightforward to check that $V({\mbox{\rv}}, \rho, M)$ satisfies
the statements of Lemma \ref{1.3.1} and Corollary \ref{1.3.2}. By
a change of variable we obtain $V({\mbox{\rv}}, \rho,
M)=\,\frac{M}{M_{\cal B}}\int_{{\cal B}}
\frac{\lambda(u)}{\|{\mbox{\rv}}-\rho u\|}\,du$.

 The results in section 2.2 can be generalized as follows.
Let ${\cal{C}}_{\rho}$ be the fixed homogeneous circle on the
$xy$-plane, centered at the origin, with radius $\rho$.
 For $\rho = 1$, write ${\cal{C}}_{1}={\cal{C}}$  and we consider $\cal C$
with fixed mass $M_{0}$.  Let $\cal B$ be such that:

(1) $dist({\cal{C}},{\cal{B}})=d\, >\, 0$,

(2) $({\cal{B}},\lambda )$ is symmetric with respect to the $xz$-plane and the $xy$-plane,

(3)  $0\, <\, M_{\cal B}\, <\, \infty$, where $M_{\cal B}$ is the mass of $\cal B$.
\vspace{.1in}

Let $V_{\cal{C}}({\mbox{\rv}}, \rho, M)$ be the gravitational
potential induced by ${\cal{C}}_{\rho}$ with total mass $M$ and
$V_{\cal{B}}({\mbox{\rv}}, \rho, M)=
\frac{M}{M_{\cal B}}\int_{{\cal{B}}_{\rho}}
\frac{\lambda_{\rho} (u)}{\|{\mbox{\rv}}- u\|}\, du=
\frac{M}{M_{\cal B}}\int_{{\cal B}} \frac{\lambda
(u)}{\|{\mbox{\rv}}-\rho u\|}\, du$, the potential induced by ${\cal B}_\rho$ with mass $M$.
Note that the $xz$-plane is an invariant subspace of both
 $V_{\cal{C}}({\mbox{\rv}}, \rho, M)$ and $V_{\cal{B}}({\mbox{\rv}}, \rho, M)$.
Let $V({\mbox{\rv}})=V_{\cal{C}}({\mbox{\rv}},1, M_{0})+
V_{\cal{B}}({\mbox{\rv}},1, M_{\cal B})$
be the gravitational potential induced by ${\cal{C}}\cup{\cal{B}}$.
Let also $V({\mbox{\rv}},\epsilon
)=V_{\cal{C}}({\mbox{\rv}},\frac{1}{\epsilon} ,\frac{1}{\epsilon}
M_{0})+ V_{\cal{B}}({\mbox{\rv}},\frac{1}{\epsilon},
\frac{1}{\epsilon}M_{\cal B})$ and
$W({\mbox{\rv}},\epsilon )=W_{\cal{C}}({\mbox{\rv}},{\epsilon})+
W_{\cal{B}}({\mbox{\rv}},{\epsilon})$, which is obtained from
$V({\mbox{\rv}},\epsilon)$ by translating the origin to $(\frac{1}{\epsilon},0,0)$.
We have $\nabla W({\mbox{\rv}},\epsilon )=\nabla W_{\cal{C}}({\mbox{\rv}},\epsilon )+
\nabla W_{\cal{B}}({\mbox{\rv}}, \epsilon )$.
Note that

{\small $$\nabla W_{\cal B} ({\mbox{\rv}},\epsilon )=\epsilon\int_{\cal{B}}
\frac{\epsilon {\mbox{\rv}}-u-  e_{1} }{\| \epsilon{\mbox{\rv}}-u- e_{1} \|^{3}}\,\lambda\, du.$$}

\noindent Since $dist({\cal{C}},{\cal{B}})=d\, >\, 0$, and $e_{1}=(1,0,0)\in\cal C$,
we have that $\| u-e_{1}\|\geq d$, for all $u\in\cal B$. Hence $\lim_{\epsilon\rightarrow
0}\nabla W_{\cal{B}}({\mbox{\rv}}, \epsilon )=0$ and it can be shown that this limit is
uniform on compacts. Therefore we can apply the methods used in section 2.2 to
prove that there are periodic orbits in the $xz$-plane, close to $\cal C$.
In the pictures below (from left to right), (1) $\cal B$ is a ball, (2) $\cal B$ is a
 three dimensional set with rotational symmetry, (3) $\cal B$ can be chosen
to be any two of the three homogeneous circles, (4) $\cal B$ is an annulus in the $xy$-plane.
\begin{figure}[!htb]
 \centering
\includegraphics[scale=0.36]{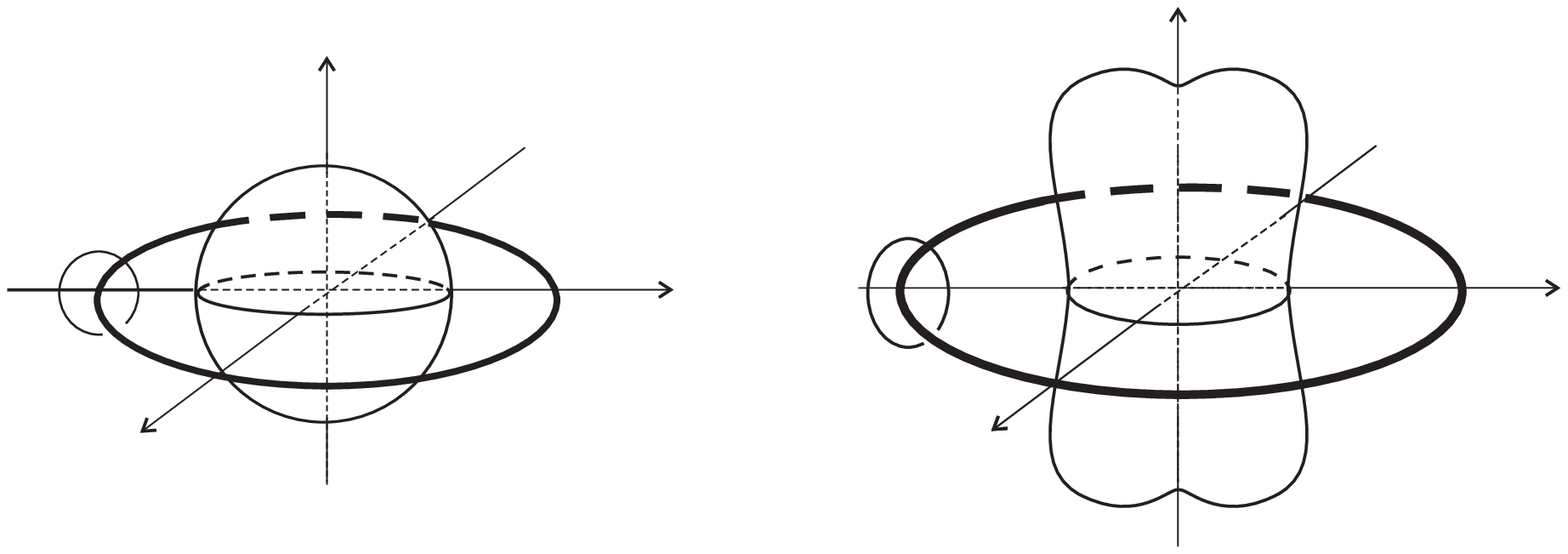}
\hspace{0.5cm}
\includegraphics[scale=0.36]{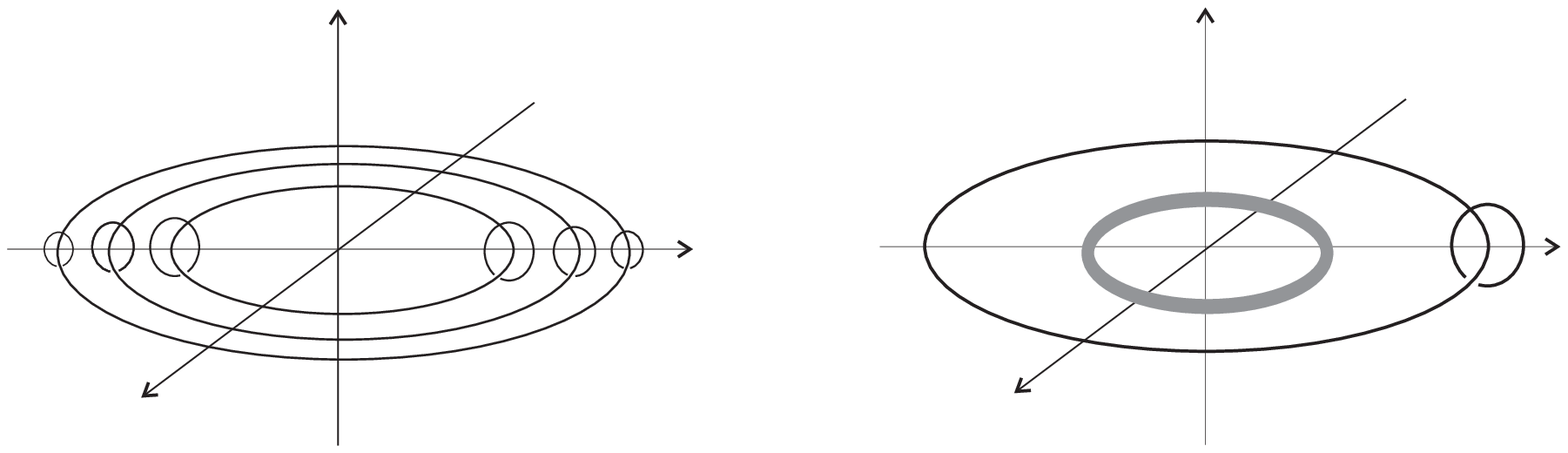}
 \end{figure}

 Analogously, it is not difficult to show that the results in section 3 can
 also be generalized. For exemple we can consider  $\cal C$ as above and
(1) $\cal B$  is the union of two inner circles and
$\cal C$ is the outer circle,
(2) $\cal B$ is annulus.

\begin{figure}[!htb]
 \centering
 \includegraphics[scale=0.35]{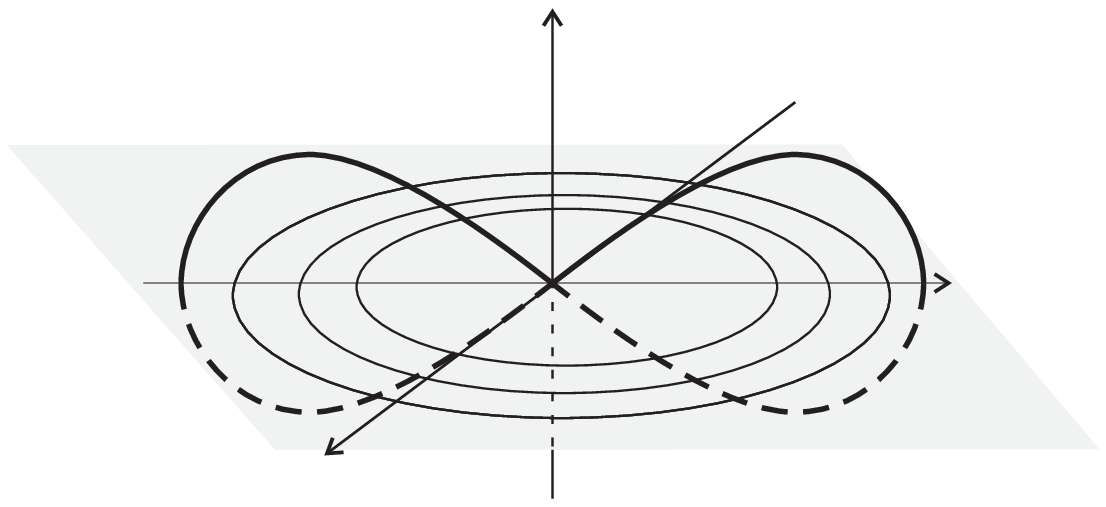}
 \hspace{0.5cm}
  \includegraphics[scale=0.35]{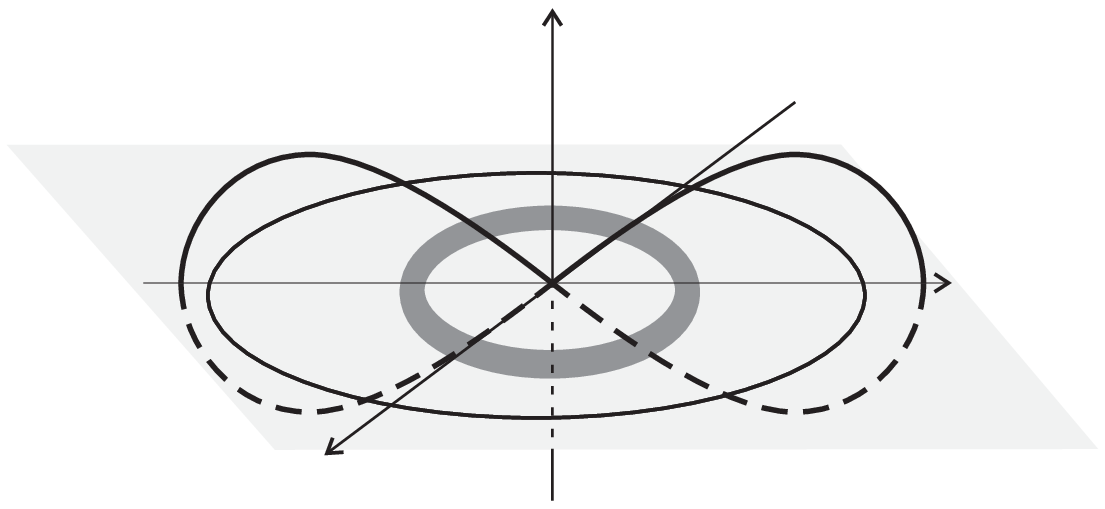}
 \end{figure}

  Proceeding as in section 3 we can prove the existence of infinitely many figure
eight periodic solutions for these two cases.

Finally, repeating the process used in section 5 we can prove the
existence of spiral solutions for bodies ${\cal C}\cup \cal B$, where $\cal B$ has
 rotational symmetry
around the $z$-axis and satisfies $dist({\cal{C}},{\cal{B}})=d\, >\, 0$.
 Here are some examples.
\begin{figure}[!htb]
 \centering
\includegraphics[scale=0.35]{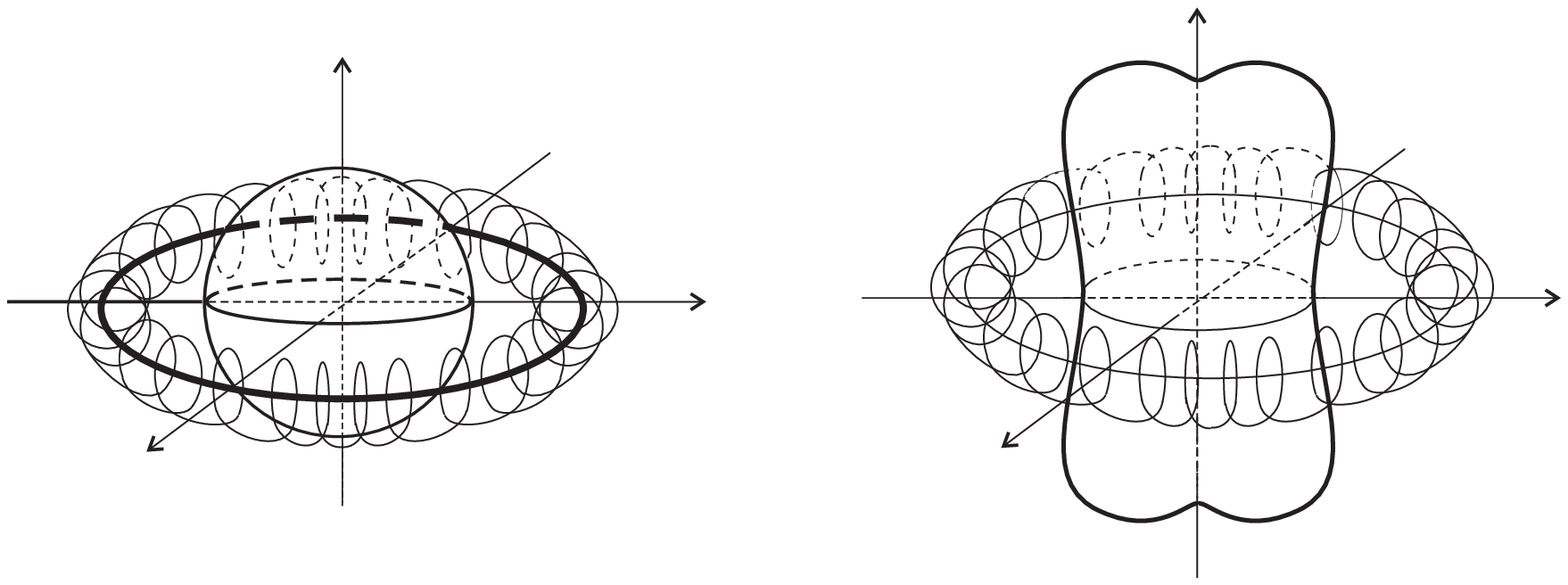}
\hspace{0.5cm}
\includegraphics[scale=0.35]{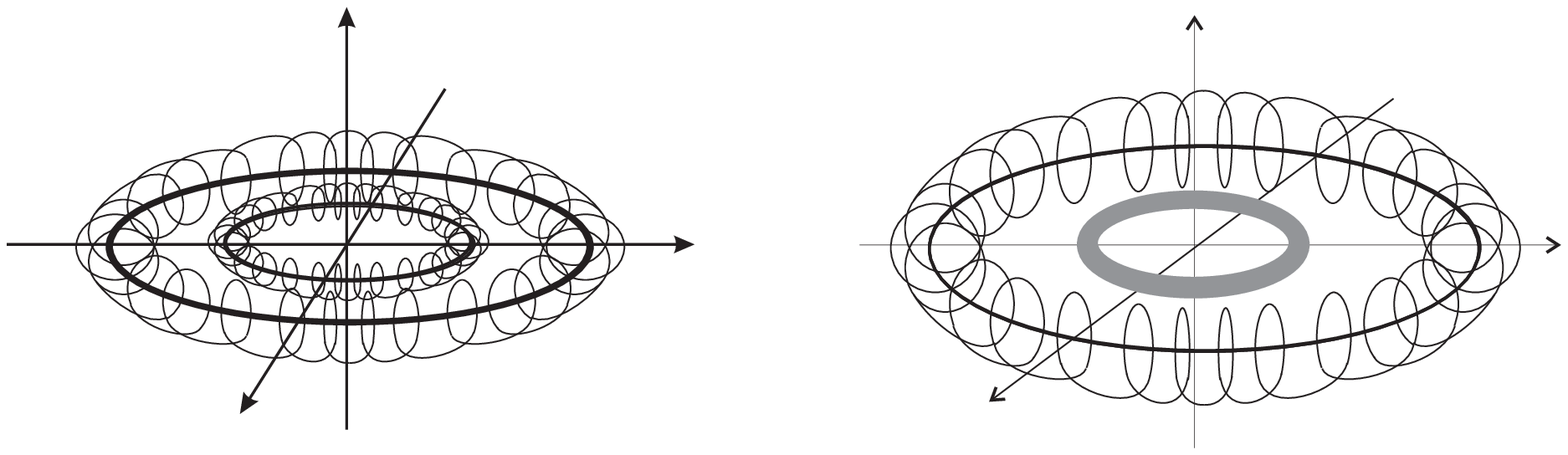}
 \end{figure}

\appendix

\section{Interval Pointing Forces.}

We identify $\R$ with $\R\times \{\, 0 \}\subset
\{(x_1, x_2 ),\,x_2\geq 0\}\subset\R^{2}$. We will consider open sets $\Omega
\subset\R^{2}$ satisfying $\hR\subset\Omega $.
Let $I\subset \R\subset\R^{2}$ be an interval and
${\mbox{\xx}}\in \chR$, ${\mbox{\vv}}\in\R^{2}$. We say that
${\mbox{\vv}}$  points to $I$ at ${\mbox{\xx}}$ if either
${\mbox{\vv}} =0$ or the infinite ray
$R({\mbox{\xx}},{\mbox{\vv}})$ that begins at ${\mbox{\xx}}$ and
has direction ${\mbox{\vv}}$, intersects $I$, i.e $
R({\mbox{\xx}},{\mbox{\vv}} )\cap I\neq \emptyset$.
{\obss {\rm \label{a}
(1) Note that if ${\mbox{\xx}}\in \R -I$ and ${\mbox{\vv}}$
points to $I$ at ${\mbox{\xx}}$ then ${\mbox{\vv}}$ is
horizontal.

\noindent(2) The following special case will be used later. Suppose
$I=(-\infty , 0]$. If ${\mbox{\vv}} =(v_{1},v_{2})\neq 0$ and
${\mbox{\xx}}=(x_{1},x_{2})\in\hR$ then ${\mbox{\vv}}$  points to
$I$ at ${\mbox{\xx}}$ if and only if
(i) $\langle {\mbox{\vv}}\, ,\, {\mbox{\xx}}^{\bot
}\rangle  =\, x_{1}v_{2}-x_{2}v_{1}\geq 0$, (ii) $v_{2}<0$.}}
 Here $(a,b)^{\bot}=(-b,a)$.
\vspace{.1in}

Statement (i) says
that the (oriented) angle from ${\mbox{\vv}}$ to $-{\mbox{\xx}}$
is non-negative, and statement (ii) says that ${\mbox{\vv}}$
points downward.
 Let $\alpha (t)$ be a curve in
$\hR$. We say that $\alpha $ points to $I$ at $t=t_{0}$ if $\dot
\alpha (t)$ points to $I$ at $\alpha (t)$. We say that $\alpha $
points to $I$ at $t\in J$ if $\alpha $ points to $I$ at $t$, for
all $t\in J$. We say that  $\alpha $ points to $I$ if $\alpha $
points to $I$ at $t$, for all $t$ in the domain of $\alpha$.

\begin{figure}[!htb]
\hspace{1cm}
\begin{minipage}[b]{0.17\linewidth}
 \includegraphics[width=\linewidth]{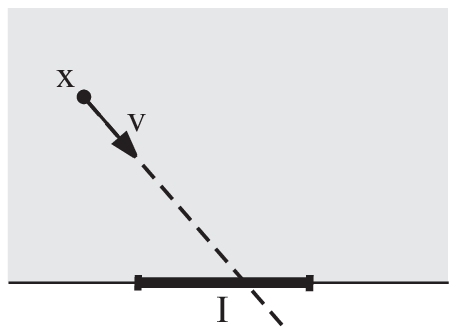}
\caption{\scriptsize{v points to I at x.}}
 \end{minipage}  \hfill
\begin{minipage}[b]{0.17\linewidth}
 \includegraphics[width=\linewidth]{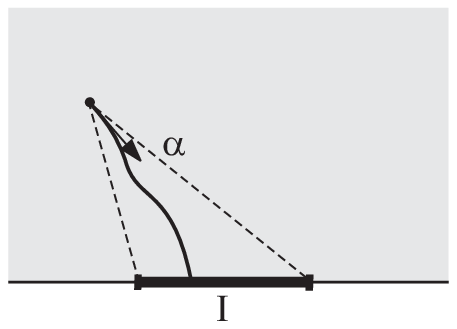}
\caption{\scriptsize{$\alpha$ points to I.}}
 \end{minipage}  \hfill
 \begin{minipage}[b]{0.17\linewidth}
 \includegraphics[width=\linewidth]{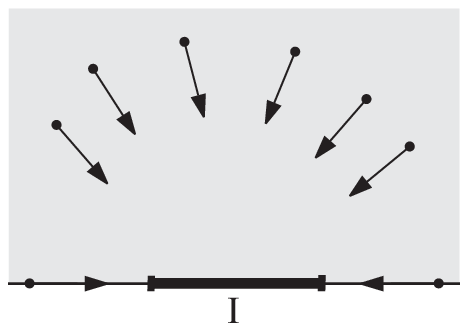}
 \caption{\scriptsize{F points to I.}}
 \end{minipage}  \hspace{1cm}
 \end{figure}

Let $\Omega$ be an open subset of $ \R^{2}$ satisfying $\hR
\subset\Omega$. Let $F:\Omega\rightarrow\R^{2}$. We say that $F$
points to $I$ if $F({\mbox{\xx}})\neq 0$ and
$F({\mbox{\xx}})$ points to $I$ at ${\mbox{\xx}}$, for every ${\mbox{\xx}}\in\hR$.



{\obss\rm{ \label{b}

(1) If $F$ is continuous and points to $I$ and ${\mbox{\xx}}\in
\R\cap\Omega$ then $F({\mbox{\xx}})$ is horizontal (but it may
happen that $F({\mbox{\xx}})=0$). Hence $\R\cap\Omega$ is an
invariant subspace of $F$.

\noindent(2) In all definitions above if an object points to
$I$ and $I\subset I' $ then it also points to $I' $.}}

\noindent(3) Let $\cal C$ be a homogeneous fixed circle in the
$xy$-plane, centered at the origin with radius $\rho$. Let $V$ be the
potential function induced by $\cal C$. It is not difficult
to see that $\nabla V$ (restricted to the $xz$-plane) points to
$[-\rho,\rho ]$.

{\prop \label{c} Let $\Omega$ be an open subset of $ \R^{2}$
with $\hR \subset\Omega$. Let $I\subset \R\subset\R^{2}$ be
a closed interval and let $F:\Omega\rightarrow\R^{2}$ be a
continuous map which points to $I$. Let also
${\mbox{\rv}}:[0,b)\rightarrow\Omega\cap \hR$ be a solution of
$\,\stackrel{..}{\mbox{\rv}}\, =\, F({\mbox{\rv}})$.
Suppose that ${\mbox{\rv}}(t)=(x(t),z(t))$ points to $I$ at $t=0$.
Then ${\mbox{\rv}}$ points to $I$.
Moreover, if ${\mbox{\rv}}$ extends to a solution
${\mbox{\rv}}:[0,b]\rightarrow\Omega\cap\chR$, with
${\mbox{\rv}}(b)\in\Omega\cap\R$, then ${\mbox{\rv}}(b)\in I$.}
\vspace{0.25cm}

\noindent {\bf Proof.} First we prove the Proposition for
$I=(-\infty , 0]$.
Since $ \stackrel{..}{\mbox{\rv}}\, =\, F({\mbox{\rv}})$ and $F$
points to $I$ by Remark \ref{a} (2) we have that
{\footnotesize \begin{equation}\left\{ \begin{array}{c} \langle
\stackrel{..}{\mbox{\rv}}(t)\, ,\, {\mbox{\rv}}^{\bot }(t)\rangle
=\,
x(t)\stackrel{..}{z}(t)-z(t)\stackrel{..}{x}(t)\geq 0\\ \\
\, \stackrel{..}{z}(t)\, <\,
0\end{array}\right.\label{p1.1}\end{equation}}

\noindent for all $t\in [0,b)$. Let $h(t)=\, \langle
\dot{{\mbox{\rv}}}(t)\, ,\, {\mbox{\rv}}^{\bot}(t)\rangle =\,
x(t)\dot{z}(t)-z(t)\dot{x}(t)$. Since
${\mbox{\rv}}(t)=(x(t),z(t))$ points to $I$ at $t=0$, we have
either $\dot{{\mbox{\rv}}}(0)=0$ or $
h(0)\geq 0,\, \dot{z}(0)\, <\, 0$.
In any case $h(0)\geq 0$. Therefore, differentiating $h$ and
using (\ref{p1.1}) we obtain that $\dot{h}(t) \geq 0$. Hence
$h(t)\geq 0$, for all $t\in [0,b)$. This proves that $ \langle
\dot{{\mbox{\rv}}}(t)\, ,\, {\mbox{\rv}}^{\bot}(t)\rangle =\,
x(t)\dot{z}(t)-z(t)\dot{x}(t)\geq 0$, for all $t\in [0,b)$.
Also, since (by \ref{p1.1}) $\stackrel{..}{z}(t)\, < \, 0$ and
$\dot{z}(0)\, <\, 0$ we have that  $\dot{z}(t)\, <\, 0$, for all
$t\in (0,b)$. Then, by remark \ref{a} (2), it follows that
${\mbox{\rv}}$ points to $I$.
This proves the first part of the Proposition (for $I=(-\infty ,
0]$). Suppose now that ${\mbox{\rv}}$ extends to a solution
${\mbox{\rv}}:[0,b]\rightarrow\Omega\cap\chR$, with
${\mbox{\rv}}(b)\in\Omega\cap\R$. Let us assume that
${\mbox{\rv}}(b)\notin I$. Since ${\mbox{\rv}}$ is continuous
$\dot{{\mbox{\rv}}}(b)$ points to $I$ at ${\mbox{\rv}}(b)$. Then,
by remark \ref{a} (1), $\dot{{\mbox{\rv}}}(b)$ is horizontal.
Since $\Omega\cap\R$ is an invariant subspace (see remark \ref{b}
(1)) it follows that ${\mbox{\rv}}(t)\in\Omega\cap\R$, for all
$t\in[0,b]$. A contradiction since ${\mbox{\rv}}(t)\in\hR$, for
$t\in [0,b)$. This proves the Proposition for the case $I=(-\infty
, 0]$. Using translations and reflections, we can prove that the
Proposition also holds for intervals $I=(-\infty , a]$, $I=[a,
\infty ) $. For the case $I=[a,b]$ apply the Proposition to
$(-\infty , b]$ and $[a, \infty ) $. This proves the Proposition.
\CaixaPreta
\vspace{0,2cm}

{\prop Let $\Omega$ be an open subset of $ \R^{2}$ with
$\hR \subset\Omega$. Let $I\subset \R\subset\R^{2}$ be a closed
interval and let $F:\Omega\rightarrow\R^{2}$ be a continuous map
which points to $I$. Let also
${\mbox{\rv}}:[0,b)\rightarrow\Omega\cap \chR$ with
${\mbox{\rv}}(0)\in\R\setminus I$ and
${\mbox{\rv}}(0,b)\subset\hR$ such that ${\mbox{\rv}}|_{(0,b)}$
is a solution of $
 \,\stackrel{..}{\mbox{\rv}}\, =\, F({\mbox{\rv}})$.
Then ${\mbox{\rv}}$ is one-to-one.}
\vspace{0.25cm}

\noindent {\bf Proof.} Note that the hypothesis:
${\mbox{\rv}}(0)\in\R\setminus I$ implies that we can assume that
$\R\neq I$. Since $F$ points to $I$ we can assume that $F$
points to a semi-infinite closed interval and without loss of
generality we can also assume that $I=(-\infty , 0]$.

 Write ${\mbox{\rv}}=(x,z)$. Since $F$ is continuous we
have that ${\mbox{\rv}}$ is a $C^{1}$ map.
Note that  $z(t) >0$ for $t\in (0,b)$.
It follows from this and the fact that $F$ points to $I$
that (see remark \ref{a} (2)) ${\mbox{\rv}}(t)=(x(t),z(t))$
satisfy \ref{p1.1} for all $t\in(0,b)$.

Now, if $\dot{z}(t)\neq 0$ for all $t\in[0,b)$ then $z$ is an
increasing function, thus one-to-one, and we have nothing to
prove. We suppose then that there is a $t_{0}\in (0,b)$ such that
$\dot{z}(t_{0})=0$. Since $ \stackrel{..}{z}(t)\, <\, 0$,
$z(t_{0})$ is a maximum value of $z$ and $t_{0}$ is unique, that
is, it is the only $t\in [0,b)$ where $\dot{z}$ vanishes. We also
have that $ \dot{z}(t) > 0$,
for $t\in [0,t_{0})$ and $\dot{z}(t)< 0$, for $t\in[t_{0},b)$.
\vspace{.1in}

\noindent {\bf Claim 1.} $\dot{x}(t_{0})\, <\, 0.$

\noindent {\bf Proof of Claim 1.}  Define
${\mbox{\sv}}:[0,t_{0}]\rightarrow\chR$ by ${\mbox{\sv}}(t) =
{\mbox{\rv}}(t_{0}-t)$.  Write
${{\mbox{\sv}}}=(\tilde{x}\, ,\,\tilde{z})$.
If $\dot{x}(t_{0})=0$ then  $\dot{{\mbox{\sv}}}(0)=0$.
 Hence
${\mbox{\sv}}$ points to $I$, at $t=0$. By Proposition \ref{c},
${\mbox{\rv}}(0)={\mbox{\sv}}(t_{0})\in I$, a contradiction.
Therefore $\dot{x}(t_{0})\neq 0$.

\begin{wrapfigure}[9]{r}{3.5cm}
 \centering
\includegraphics[width=3cm]{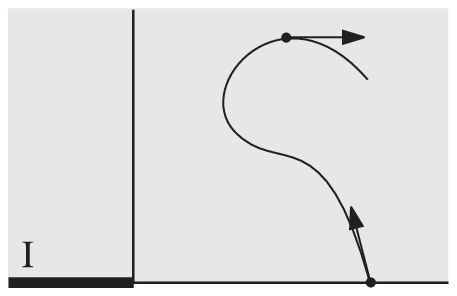}
 \caption{\scriptsize{$\dot{x}(t_{0})>0$ cannot happen.}}
 \end{wrapfigure}

Suppose that  $\dot{x}(t_{0})\, >\, 0$ (see figure A.15). Then
$\dot{\tilde{x}}(0)\, <\, 0$. Hence, for  $t$ close to $0$,
$\dot{\tilde{x}}(t)\, <\, 0$. Since $\dot{z}(t)\, >\, 0$ for $t\,
<\, t_{0}$, we have that $\dot{\tilde{z}}(t)\, <\, 0$ for $t\, >\,
0$. Therefore, since ${{\mbox{\sv}}}$ is $C^{1}$,
${{\mbox{\sv}}}$ points to $I=(-\infty , 0]$ for $t$ close
to $0$. By the Proposition above
${{\mbox{\sv}}}(t_{0})\in I$, which is a contradiction
since ${{\mbox{\sv}}}(t_{0})={\mbox{\rv}}(0)\notin I$. This
proves the claim.

By the claim above, for  $t$ close to $t_{0}$, $\dot{x}(t)\, <\,
0$. Since we have that $\dot{z}(t)\, <\, 0$ for $t\, >\, t_{0}$
and ${\mbox{\rv}}$ is $C^{1}$, it follows that ${\mbox{\rv}}$ points
to $I=(-\infty , 0]$ for $t$ close to $t_{0}$ and $t\, >\, t_{0}$.
Then, by the Proposition above, ${\mbox{\rv}}$ points to $I$ at
$t\in (t_{0}, b)$.

Let ${\mbox{\sv}}(t) = {\mbox{\rv}}(t_{0}-t)$ as in proof of the claim 1.
Suppose ${\mbox{\rv}}$ is not one-to-one. Let
$\bar{t}=min\{ \, t\in (t_{0},b)\,\,
:\,\,{\mbox{\rv}}(t)={\mbox{\rv}}(t'){\mbox{ for some
}}t'\in[0,t_{0})\,\}$.
Note that $\bar{t}$ exists and $\bar{t}\, >\, t_{0}$, since
$\dot{x}(t_{0})\, <\, 0$. Let $t'$ be such that
${\mbox{\rv}}(\bar{t})={\mbox{\rv}}(t')$. Note that $t'$ is
uniquely defined and $t'\in (0,t_{0})$ because $z$ is increasing
on $[0,t_{0}]$.
Write $d=\bar{t}-t_{0}$ and $c=t_{0}-t'$ and define
$\alpha ={\mbox{\sv}} |_{[0,c]}$ and
$\beta :[0,d]\rightarrow \Omega\cap\hR,\,\beta (t)=
{\mbox{\rv}}(t-t_{0}).$

\begin{wrapfigure}[9]{r}{3cm}
 \centering
\includegraphics[width=3cm]{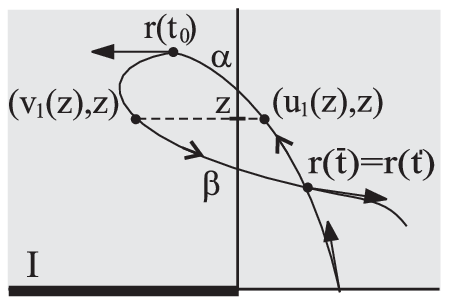}
 \caption{\scriptsize{This is what happens if ${\mbox{\rv}}$ is not one-to-one.}}
 \end{wrapfigure}

Write $\alpha =(u_{1},u_{2})$ and $\beta = (v_{1},v_{2})$. Then
$\alpha$ and $\beta$ have the following
 properties: (property (d)
follows from the minimality of $\bar{t}$ and claim 1)

(a) $\alpha $ does not point to $I$ at $t$, for all $t\in [0,c]$
and  $\beta (t)$ points to $I$ at $t\in(0, d]$,

(b) $\alpha (0)=\beta (0)$ and $\alpha (c)=\beta (d)$,

(c) $\dot{u_{2}}(t)\, <\, 0$ for $t\in(0,c]$ and $\dot{v_{2}}(t)\,
<\, 0$ for $t\in(0,d]$,

(d) If $u_{2}(t)=v_{2}(t^{*})$ then $v_{1}(t)\, <\, u_{1}(t^{*})$,
for $t\neq 0,c$.

\vspace{0.25cm}

\noindent {\bf Claim 2.}{\it  The (oriented) angle from
$\dot{\alpha}(c)$ to $\dot{\beta}(d)$ is non-negative and less
than $\pi$.}

\noindent {\bf Proof of Claim 2.} First, since both vectors point
downward, we have that this angle is less that
$\pi$.
Write $z_{0}= u_{2}(0)=v_{2}(0)=z(t_{0})$, and $z_{1}=
u_{2}(c)=v_{2}(d)=z(t')=z(\bar{t})$.
By property (c) above we can write $u_{1}$ and $v_{1}$ in terms
of $z\in[z_{1},z_{0})$. That is, there is $t=t(z)$, $z\in[z_1, z_0)$
such that  $u_{2}(t(z))=z$. Then we write $u_1(z)=u_1(t(z))$.
Similarly we can write $v_1(z)$.
Note that, by properties (b) and (d),
$u_{1}(z_{1})-v_{1}(z_{1})=0$, and  $u_{1}(z)-v_{1}(z)\geq 0$, for
$z\in(z_{1},z_{0})$. Hence,
$\frac{d}{dz}(u_{1}(z)-v_{1}(z))\mid
_{z_{1}}\geq 0$.
 But $\frac{d}{dz}(u_{1}(z))\mid
_{z_{1}}=\frac{\dot{u_{1}}(c)}{\dot{u_{2}}(c)}$ and
$\frac{d}{dz}(v_{1}(z))\mid
_{z_{1}}=\frac{\dot{v_{1}}(d)}{\dot{v_{2}}(d)}.$ Therefore
$\frac{\dot{u_{1}}(c)}{\dot{u_{2}}(c)}\geq
\frac{\dot{v_{1}}(d)}{\dot{v_{2}}(d)}$ and it follows that
$ \langle \dot{\beta}(d)\, ,\, \dot{\alpha}^{\bot}(c)  \rangle \geq 0$.
This proves claim 2. \vspace{.1in}

By the claim above and property (a),
$\alpha$ does not point to $I$ at $c$. It follows that $\beta$ does
not point to $I$ at $d$. This contradicts property (a) and proves
the Proposition.
\CaixaPreta
\vspace{0.6cm}

\noindent {\bf Remark.} Let us assume that $I=[a,b]\,$ is bounded.
If we shoot a particle $P$ upward from the $x$-axis, by the
Proposition above (the connected piece in $\chR$ of) the solution
behaves, with respect to self-intersections and landing, in the
following way.
If we shoot from outside $I$, say, to the right of $I$, we have no
self-intersections.
\begin{figure}[!htb]
 \centering
\includegraphics[scale=0.6]{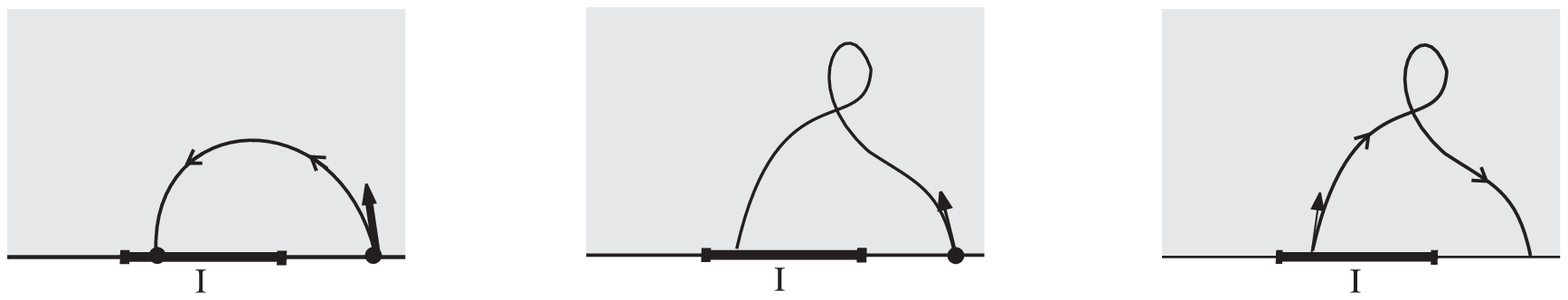}
\caption{\scriptsize{The figure to the left can happen. The two figures
to the right cannot happen.}}
 \end{figure}

Also, the particle can land on $I$, or to the left of $I$, but not
to the right of $I$. (To see that the particle $P$ cannot land to
the right of $I$ see the proof of claim 1 in the proof of the
Proposition above, and use also the fact that the orbit has no
self-intersections.)
\begin{figure}[!htb]
 \centering
\includegraphics[scale=0.6]{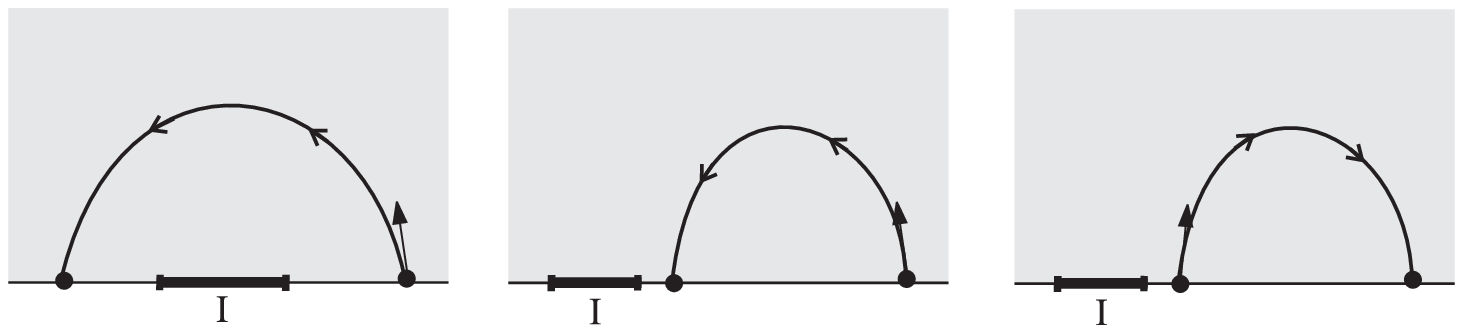}
 \caption{\scriptsize{The figure to the left can happen. The figure in the middle
  cannot happen. The figure to the right cannot happen.}}
 \end{figure}

If we shoot from inside $I$ then we may have self-intersections
but then the particle should land inside $I$, otherwise, inverting
time, we obtain one of the previous cases. Also, in this case, if
the orbit of $P$ has no self-intersections then it may land
anywhere.
Of course, the particle may not land on the $x$-axis at all, and
escape to infinity.
The same analysis can be made for closed semi-infinite intervals
$I$.
\vspace{.2cm}

\section{Proof of Proposition \ref{3.8}.}

To  prove Proposition \ref{3.8} we will use the  following Lemmas
and remarks. Write  $\Lambda_A
= \frac{2}{min\{1,A\}} +1,$ for $A>0.$ Note that $\Lambda_A >1.$
In the following two Lemmas we assume $z(t)$ to be twice differentiable.

{\lem Let $z(t)$ be defined in $[0,a]$  and such that
(i) $\stackrel{..} z \,\,\leq
-Az,\,A>0$, (ii) $z(0)=z_{0}>0$, (iii) $\dot{z}(0)=0.$
If $a>\Lambda_A$ then there is $t_{0}\in [0, \Lambda_{A}]$  such
that $z(t_{0})=0$. \label{0.1.2}}

\vspace{0,2cm}

\noindent{\bf Proof.} Suppose that $z(t)\neq 0$ for all $t\in [0,
\Lambda_{A}]$. Then   $z(t)>0$ and  $\dot{z}(t)$  is decreasing  in
$ [0, \Lambda_{A}]$.
Since $\dot{z}(0)=0,$ we have $\dot{z}(t)<0,$ $\,t\in
(0,\Lambda_{A}]$. Then    $z$ is decreasing  in $(0,\Lambda_{A}].$
Note  also that $z$ is defined in $1$ because $1<\Lambda_A <a$.
We claim that $\dot{z}(1)\leq -\bigl( {min\{1,A\}}
\bigr)\frac{z_{0}}{2}.$ We have two  cases: $z(1)\geq
\frac{z_{0}}{2}$ or $z(1)\leq \frac{z_{0}}{2}$.
Suppose first that $z(1)\geq \frac{z_{0}}{2}$. This implies that
$z(t)\geq \frac{z_{0}}{2}$ for all $t\in [0,1]$. Then

{\footnotesize $$\dot{z} (1)  =  \int_{0}^{1} \stackrel{..} z (t) dt \leq
 \int_{0}^{1} -A z (t)dt
 \leq  -A \int_{0}^{1}  z (t)dt \leq - A \frac{z_{0}}{2}.$$}

Suppose now that $z(1)\leq \frac{z_{0}}{2}$. By the intermediate
value theorem there exists  $t' \in(0,1)$ such  that
 {\small $\dot{z}(1) < \dot{z}(t') = \frac{z(1)-z(0)}{1-0}
 \leq -\frac{z_{0}}{2}.$}
This proves our claim.

Since $\dot{z}$ is decreasing  in $(0,\Lambda_{A}]$, for $1\leq
t\leq \Lambda_{A}$ we have that   $\dot{z}(t)\leq \dot{z}(1)\leq
-\alpha$, with $\alpha= (min\,\{A,1\})\frac{z_{0}}{2}.$
Hence, $z(t)-z(1)=\int_{1}^{t} \dot{z}(s) \,ds\leq
-\alpha\,(t-1)$. It follows  that $z(t)\leq z(1) -\alpha\,(t-1)\leq
z_{0}-\alpha \,(t-1)$, for $1\leq t\leq \Lambda_{A}$.  Therefore
{\small $z(\Lambda_{A})=z\Bigl( \frac{z_{0}}{\alpha} +1\Bigr) \leq z_{0}-\alpha
\,(\Lambda_A -1)=0$}, which is a contradiction.
 It follows  that there exists  $t_0 \in
[0,\Lambda_A ]$ such  that $z(t_0)=0$. \CaixaPreta

\vspace{0,3cm}

{\lem Let $z(t)$ be  defined  in $[0,+\infty)$ and such  that
(i) $\stackrel{..}z \,\,\leq -Az$, $A>0$, (ii) $z(0)=0$,
 (iii) $\dot{z}(0)>0.$
Then  there exists  $t_{0}>0$ such  that $\dot{z}(t_{0})=0$.
\label{0.1.3}}

\vspace{0,2cm}

\noindent{\bf Proof.} Suppose that      $\dot{z}(t)\neq 0$, for
all $t$. Then    $\dot{z}(t)>0$ for all $t$. Hence $z$  is an
increasing function.
Then for all $t\geq 1$, $0<z(1)\leq z(t)$  and
{\footnotesize $\dot{z}(t)=\dot{z}(1)+\int_{1}^{t} \stackrel{..} z (s)ds \leq
\dot{z}(1)-A \int_{1}^{t}  z (s)ds \leq \dot{z}(1)-A z (1)(t-1).$}
 Evaluating at $t= 1+\frac{\dot{z}(1)}{A {z}(1)}$ we have
$\dot{z}(1+\frac{\dot{z}(1)}{A {z}(1)})\leq 0,$ which is a
contradiction. \CaixaPreta

\vspace{0,2cm}

Before proving  Proposition \ref{3.8} we have the  following
remark.

 {\obss {\rm  \label{c3}

\noindent   Let $\delta<0.$ By Lemma 2.3 of \cite{AO1}, there
exists  $R_{\delta}$, $0<R_{\delta}<+\infty$, such  that
$\{\,\mbox{\rv}\,;\,V(\mbox{\rv})\leq \delta\,\}\subset
B(R_{\delta})$  where $B(R_{\delta})$ is the ball centered at the
origin of radius $R_{\delta}$. Therefore, if $\mbox{\rv}(t)$ is a
solution of $\,\stackrel{..} {\mbox{\rv}}\,=-\nabla
{V}(\mbox{\rv})$
 (restricted to the
$xz$-plane) with $E(\mbox{\rv}(t))\leq \delta<0,$
then $\|\mbox{\rv}(t)\|\leq R_{\delta}.$ Moreover, if
  $\cal C$ is the fixed homogeneous
circle of radius one and centered at the origin and $u\in\cal C$, we have:
$\|\mbox{\rv}-u\|\leq \|\mbox{\rv}\|+ \|u\|\leq R_{\delta}+1$.
 Hence,
$\frac{1}{\|\mbox{\rv}-u\|^{3}}\geq \frac{1}{(
R_{\delta}+1)^{3}}$, and follows that
{ $
\lambda   \int_{\cal C}\frac{du}{\| \mbox{\rv}-u\|^{3}}\geq
\lambda \int_{\cal C}\frac{du}{ ( R_{\delta}+1)^{3}}  = \frac{M}{
( R_{\delta}+1)^{3}}.$}}}

\vspace{0,4cm}

\noindent{\bf Proof of Proposition \ref{3.8}.} From (\ref{8.1.4}) we have
{\small
  $\stackrel{..} z \,\,=-\frac{\partial V}{\partial z} (x,z)=-\lambda
  z\int_{\cal C}\frac{du}{\| \mbox{\rv}-u\|^{3}}$}.
Note that $ \stackrel{..} z\,\, <0,$ if $z>0$ and $ \stackrel{..}
z \,\,>0$, if $z<0$. Let $\mbox{\rv}(t)=(x(t),z(t))$ be a solution of
$\,\stackrel{..} {\mbox{\rv}}\,=-\nabla {V}(\mbox{\rv})$ (restricted to the
$xz$-plane),
 with $E(\mbox{\rv}(t))\leq \delta<0,$ $z(0)=0$ and
$\dot{z}(0)>0$. Then $| \stackrel{..} z |=\lambda \,|z| \int_{\cal
C}\frac{du}{\| \mbox{\rv}-u\|^{3} }\geq
\frac{M}{(R_{\delta}+1)^{3}}\,|z|$ (see Remark \ref{c3}). Set
$A=\frac{M}{(R_{\delta}+1)^{3}}$ and $\Lambda_A =
\frac{2}{min\{1,A\}} +1.$ We have $| \stackrel{..} z |\geq
A\,|z|$. If $z>0$, then $ \stackrel{..} z\,\,<0$, and the
inequality becomes $ -\stackrel{..} z \,\,\geq Az$, for $z\geq 0$.
Let $(a,b),$ $a<0<b$, be the maximal interval on which
$\mbox{\rv}(t)$ is defined. If $b<\infty$, $\mbox{\rv}(t)$
collides (see \cite{AO1}, Theorem A), hence $\lim_{t\rightarrow
b^-} z(t)=0,$ and we have nothing to prove. Suppose then that
$b=\infty$. By Lemma \ref{0.1.3} there exists  $\bar t \in(0,\infty),$
such  that $\dot{z}(\bar t)=0.$
Replacing $\bar t$ by $min
\{\,t>0\,;\,\dot{z}(t)=0\,\}$, we can assume that $\dot z (t)>0,$
$\,t\in[0,\bar t\,).$ Then $z(t)$ is increasing in $[0,\bar t\,] $
and $z(t)>0,$ $\,t\in(0,\bar t\,]$. Applying Lemma \ref{0.1.2} to
the function $t\mapsto z(\bar t -t),$ $t\in[0,\bar t],$ we obtain
that  $\bar t \leq \Lambda_A$.
 We will show
that there exists  $t_0 \in [0,2\Lambda_a ],$ such  that
$z(t_0)=0$. Suppose that $z(t)\neq 0$, $t\in (0,2\Lambda_A ]$.
Since $\bar t \in (0, 2\Lambda_A )$ and $z(\bar t)>0$ we have
$z(t)>0$, $t\in(0,2\Lambda_A].$ Hence there exists  $c$ such
that $z(t)>0,$ $\,t\in(0,c),$ with $2\Lambda_A <c<\infty$. Note that
$c-\bar t >\Lambda_A.$
Applying Lemma \ref{0.1.2} to the function $t\mapsto z(t+\bar t),$
$\,t\in [0,c-\bar t]$ we obtain a contradiction. This  proves the
Proposition. \CaixaPreta

\vspace{0,35cm}

{}

\vspace{.5in}

P. Ontaneda and C. Azev\^edo,

e-mail: ontaneda@dmat.ufpe.br

Departamento de Matem\'atica,

Universidade Federal de Pernambuco

Recife, PE, 50670-901

Brazil

\end{document}